%% file: paper_rquf.tex
\RequirePackage[2018-12-01]{latexrelease}
\documentclass{rQUF2e}
\usepackage{hyperref}
\makeatletter
\let\orig@ref\ref
\renewcommand{\ref}[1]{%
  \expandafter\@setref\csname r@#1\endcsname\@firstoffive{#1}%
}
\def\@firstoffive#1#2#3#4#5{#1}
\AtBeginDocument{%
  \let\hyperref@ref\ref
  \renewcommand{\ref}[1]{%
    \expandafter\@setref\csname r@#1\endcsname\@firstoffive{#1}%
  }
}
\makeatother
\usepackage{natbib}
\usepackage{float}
\usepackage{caption}
\usepackage{graphicx}
\usepackage{epstopdf}
\usepackage{subfigure}
\usepackage{subcaption}
\usepackage{tikz}
\usetikzlibrary{patterns,shadows.blur}
\usepackage{url}

\usepackage{xcolor}

\newcommand{\layer}{\ell}
\newcommand{\loss}{\mathcal{L}}
\newcommand{\attach}{A}
\newcommand{\limit}{L}
\newcommand{\shift}{\Delta}
\newcommand{\reinst}{r}
\newcommand{\premium}{\Pi}
\newcommand{\recovery}{R}

\newcommand{\peril}{p}
\newcommand{\groups}{\mathcal{G}}
\newcommand{\group}{g}

\newcommand{\subgroup}{s}
\newcommand{\trial}{t}
\newcommand{\event}{e}
\newcommand{\rol}{\mathrm{rol}}
\newcommand{\lol}{\mathrm{lol}}
\newcommand{\profit}{\pi}
\newcommand{\objective}{\mathcal{O}}

\newcommand{\netprofit}[1]
{\profit_{\mathrm{net},#1}}
\newcommand{\grossprofit}{\Pi_{\text{gross}}}
\newcommand{\netloss}{\mathcal{N}}

\newbox{\orcidicon}
\sbox{\orcidicon}{\raisebox{-0.5ex}{\includegraphics[height=2.2ex]{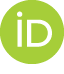}}}
\newcommand{\orcid}[1]{%
  \begingroup%
  \hypersetup{pdfborder={0 0 0}}%
  \href{https://orcid.org/#1}{%
    \kern-1pt\usebox{\orcidicon}\kern1pt%
  }%
  \endgroup%
}

\begin{document}
%\jvol{00} \jnum{00} \jyear{2014} \jmonth{October}

\title{
  Modern Computational Methods in Reinsurance Optimization: From Simulated
  Annealing to Quantum Branch \& Bound\\
}

\author{
  George~Woodman$^{\ast}\,$\orcid{0009-0002-3403-4201}$^{1}$,
  Ruben~S.~Andrist$^{\ast}\,$\orcid{0000-0003-1126-6931}$^{2}$,
  Thomas~H\"aner~\orcid{0000-0002-4297-7878}$^{2}$,
  Damian~S.~Steiger~\orcid{0000-0003-1588-8930}$^{2}$,
  Martin~J.~A.~Schuetz~\orcid{0000-0001-5948-6859}$^{2}$,
  Helmut~G.~Katzgraber~\orcid{0000-0003-3341-9943}$^{2}$,
  Marcin~Detyniecki~\orcid{0000-0001-5669-4871}$^{1}$
  \thanks{
  $^\ast$Corresponding authors: george.woodman@axa.com, randrist@amazon.com}
}

\affil{
  \label{affiliation:axa}$^{1}$AXA Group Operations, Paris, France\\
  \label{affiliation:qsl}$^{2}$Amazon Advanced Solutions Lab, Seattle, Washington 98170, USA
}
\received{Dated: \today}

\maketitle

\begin{abstract}
  We propose and implement modern computational methods to enhance catastrophe
  excess-of-loss reinsurance contracts in practice. The underlying optimization
  problem involves attachment points, limits, and reinstatement clauses, and
  the objective is to maximize the expected profit while considering risk
  measures and regulatory constraints. We study the problem formulation,
  paving the way for practitioners, for two very different approaches: A local
  search optimizer using simulated annealing, which handles realistic
  constraints, and a branch \& bound approach exploring the potential of a
  future speedup via quantum branch \& bound.

  On the one hand, local search effectively generates contract structures
  within several constraints, proving useful for complex treaties that have
  multiple local optima. On the other hand, although our branch \& bound
  formulation only confirms that solving the full problem with a future quantum
  computer would require a stronger, less expensive bound and substantial
  hardware improvements, we believe that the designed application-specific
  bound is sufficiently strong to serve as a basis for further works.
  Concisely, we provide insurance practitioners with a robust numerical
  framework for contract optimization that handles realistic constraints today,
  as well as an outlook and initial steps towards an approach which could
  leverage quantum computers in the future.
\end{abstract}


\begin{classcode}
  G22;   % Insurance, Insurance Companies, Actuarial Studies
  % C15; % Statistical Simulation Methods
  C61;   % Optimization Techniques, Programming Models, Dynamic Analysis
  C63    % Computational Techniques, Simulation Methods
\end{classcode}

%%%%%%%%%%%%%%%%%%%%%%%%%%%%%%%%%%%%%%%%%%%%%%%%%%%%%%%%%%%%%%%%%%%%%%%%%%%%%%%
%%%%%%%%%%%%%%%%%%%%%%%%%%%%%%%%%%%%%%%%%%%%%%%%%%%%%%%%%%%%%%%%%%%%%%%%%%%%%%%
\section{Introduction}

% base intro
Reinsurance is a way for insurance companies to reduce their risk exposure: by
ceding parts of the losses to a third party for a premium.
In the context of catastrophe excess-of-loss reinsurance, crafting a suitable
reinsurance contract requires balancing multiple competing objectives:
maximizing expected profit, minimizing earnings volatility, adhering to
internal business standards, meeting regulatory requirements, and ensuring the
contract can be placed in current reinsurance market conditions.
This multi-dimensional optimization problem is further complicated by the
non-linear nature of reinsurance contracts, where layered structures,
reinstatement clauses, and aggregate features create complex inter-dependencies
among parameters.

% context (previous work)
The optimal solution depends on the objective function, risk measures and
pricing criteria considered. For pricing according to the expected value
principle, \cite{Borch1960} showed that stop-loss is the ideal choice, with
\cite{Arrow1963} obtaining similar results for maximizing the utility of the
insurer's final wealth. While this work was extended to several other pricing
approaches (see~\cite{Albrecher2017} for an overview), more recent publications
have focused on risk measures such as the value at risk (VaR), conditional
value at risk (CVaR), Tail value at risk (TVaR) and conditional tail
expectation (CTE). CVaR, TVaR and CTE are mathematically closely related and
are often used interchangeably within the literature.
See, e.g., \cite{Dhaene2006}, \cite{Cai2007}, and \cite{Chi2011}. For a more
complete overview, we refer to \cite{Cai2020}.

% Focus on XL + practical application
The analytical frameworks employed in these studies rely on strong assumptions
about loss distributions, limiting their practical applicability. In practice,
the conditions for obtaining reinsurance are dictated by the market and data
sets for risk evaluation are prescribed, in part, by regulatory requirements.
In this work, we focus on catastrophe excess-of-loss reinsurance contracts with
a set of features modeled after placements in recent years. Modern catastrophe
programs often involve multiple layers with varying attachment points, shared
limits, and reinstatement clauses. The interaction between these parameters
creates a complex optimization landscape with local optima, making traditional
gradient-based methods unsuitable. Furthermore, practical constraints such as
minimum retention requirements, attachment probability, and thresholds for
loss-distribution derived observables further complicate the optimization
process.

% MCMC
To bridge this gap, we propose to leverage Markov Chain Monte Carlo (MCMC)
methods~\citep{Metropolis1953}, which are widely used to numerically explore
complex optimization landscapes in practice. Simulated annealing is an MCMC
method which emerged from the connection between statistical mechanics and
combinatorial optimization. \cite{Kirkpatrick1983} showed how the approach can
be leveraged to solve difficult combinatorial problems in the area of circuit
placement for computer chips. Its versatility has led to applications across
diverse domains, ranging from spin-glass systems \citep{Marinari1992}, to
weight optimization in neural network training~\citep{Sexton1999}, vehicle
routing~\citep{Osman1993} and job shop scheduling~\citep{VanLaarhoven1992}, to
name a few.

% Quantum
In addition to using classical computers for reinsurance optimization, we also
discuss the potential of quantum computing for solving such problems with a
focus on quantum branch and bound (QBB) \citep{montanaro2020quantum}. We
recognize, that while the near-quadratic speedup of QBB over classical branch
and bound (B\&B) may seem promising, the clock speed of future fault-tolerant
quantum computers is expected to be orders of magnitude lower than current
classical clock speeds~\citep{Babbush2021QuadraticSpeedup} of CMOS hardware,
which severely restricts the range of B\&B methods that could benefit from QBB
in practice. But, as pointed out in previous work~\citep{haner2024solving},
early quantum advantage of QBB is most likely for B\&B methods featuring
computationally inexpensive (but effective) bounding operators and large trees.
This serves as our motivation for deriving such a bound and implementing an
application-specific B\&B solver for a reinsurance optimization problem.

% Contents
This paper makes three main contributions to the literature. First, we develop
the mathematical formulation and a flexible computational framework for
reinsurance optimization that bridges the gap between theoretical models and
practical implementation for two distinctive and complementary approaches.
Second, we introduce and analyze a Markov process and simulated annealing
implementation tailored for optimizing reinsurance contracts. Third, we
propose an application-specific branch \& bound approach to a reinsurance
optimization problem for a given set of peril groups, and show that our
proposed iterative bound is sufficiently strong to effectively prune a large
fraction of nodes from the branch \& bound tree.

% Organization
The remainder of this paper is organized as follows: The reinsurance model and
features considered in this paper are defined in section~\ref{sec:model}.
Section~\ref{sec:markov-chain} presents our optimization framework, details the
Markov-chain-based optimization approach and its results. The branch \& bound
approach is discussed in section~\ref{sec:branch-bound}.
Section~\ref{sec:discussion} concludes with a discussion of implications for
practical use and suggestions for future research.

%%%%%%%%%%%%%%%%%%%%%%%%%%%%%%%%%%%%%%%%%%%%%%%%%%%%%%%%%%%%%%%%%%%%%%%%%%%%%%%
%%%%%%%%%%%%%%%%%%%%%%%%%%%%%%%%%%%%%%%%%%%%%%%%%%%%%%%%%%%%%%%%%%%%%%%%%%%%%%%
\section{Reinsurance Model} 
\label{sec:model}
There are several types of reinsurance cover, which can cover different risks
and offer different benefits to both the insurer and the reinsurer. There are
two basic categories of reinsurance: treaty and facultative. Treaties are
agreements that cover broad groups of policies, such as all of an insurer's US
property business. Facultative covers specific individual, usually high-value
or high-risk risks that wouldn't be acceptable under a treaty. There are
endless complexities that can be used to make a reinsurance treaty fit for
purpose.

We focus on catastrophe treaty excess-of-loss reinsurance, which provides
insurance companies (cedants) protection against severe natural and man-made
perils. Unlike traditional reinsurance, which might cover individual claims or
aggregate losses over a period, catastrophe excess-of-loss reinsurance
specifically addresses large-scale events that can cause significant losses
across multiple policies simultaneously. These events include natural perils
such as hurricanes, earthquakes, and floods, as well as man-made catastrophes.

A treaty excess-of-loss reinsurance contract typically consists of multiple
\textit{layers} $\layer$, each covering events of specific peril type(s), such
as Floods or Windstorms. They are defined by an \textit{attachment point}
$\attach$ (where coverage begins) and a \textit{limit} $\limit$ (maximum
recovery). When losses exceed the attachment, the reinsurer pays up to the
layer's limit, while the cedant retains losses below the attachment and
exceeding the upper threshold, $A+L$. If the losses exceed $\attach$ but do not
exceed $A+L$, then the limit is only partially exhausted and the remaining
limit may be used for subsequent losses exceeding $\attach$. If the limit is
exhausted, a reinstatement clause may allow the limit to reset, providing
additional coverage for subsequent losses within the same policy period. This
is usually in exchange for a reinstatement premium. In the following we
consider a catastrophe treaty excess-of-loss reinsurance contract structure
that accommodates multiple grouped perils, combined coverage, and reinstatement
clauses.

%%%%%%%%%%%%%%%%%%%%%%%%%%%%%%%%%%%%%%%%%%%%%%%%%%%%%%%%%%%%%%%%%%%%%%%%%%%%%%%
\vspace{-5mm}
\subsubsection*{Perils and Grouping:}
We define a hierarchical grouping structure for the perils $\peril$ covered:
\begin{itemize}
  \item Groups $\group \in \groups$ partition the perils into main coverage
    units (perils that cede against a common aggregate yearly limit for each
    shared layer).
  \vspace{2mm}
  \item Within each group $\group$, subgroups $\subgroup$ can further partition
    the perils (perils that share a common attachment level).
\end{itemize}
This hierarchical structure is motivated by typical reinsurance structures: A
single excess-of-loss layer can cover multiple perils concurrently, indicating
that their recovered amounts contribute towards a common aggregate limit.
However, because the magnitude and frequency of events can differ substantially
by peril type, it may be desirable to have different attachment points for this
shared coverage. Defining the limit at the group level and attachment at the
subgroup level allows grouping of such perils with different levels of
activity. 

We use loss data that is generated by a natural catastrophe model that
simulates thousands of trial years of peril events applied to an insurance
portfolio. As a result, each gross loss event $\event$ (generating a loss
${\loss}_{\event} > 0$ in a trial year $\trial$) is associated with a peril
${\peril}_\event$.

\begin{figure}
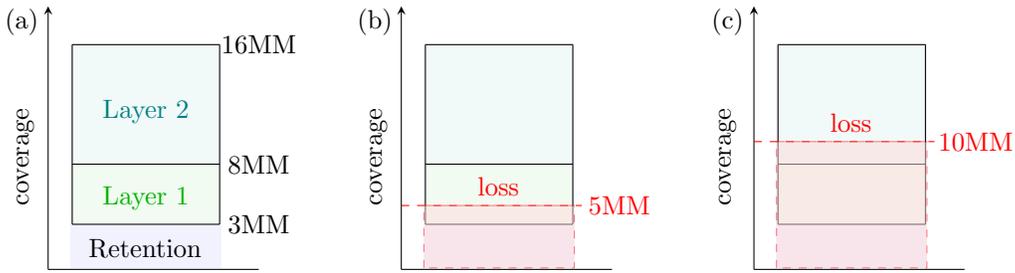

  \vspace{3mm}
  \centering
  \include{figures/layered}
  \vspace{-1mm}
  \captionsetup{width=.8\linewidth}
  \caption{
    \textbf{(a)} Layered reinsurance profile: Layer 1 has attachment 3MM and
      limit 5MM, Layer 2 has attachment 8MM and limit 8MM.
    \textbf{(b)} Distribution of 5MM from a single loss event (red): The first
      3MM in losses are retained by the cedant. The remaining 2MM are ceded to
      the reinsurer of layer 1, which is now 40\% exhausted (all of which could
      be reinstated, if designated in the contract).
    \textbf{(c)} Assuming this is the first applicable loss in the treaty
      period, the distribution of 10MM in losses between insurer and
      reinsurers: The first 3MM in losses are retained, the next 5MM are ceded
      to the reinsurer of layer 1, exhausting that layer. The remaining 2MM are
      ceded to the reinsurer of layer 2. In a typical case, if not reinstated,
      layer 1 will not cover any future claims and layer 2 still has 8MM - 2MM
      = 6MM available for future claims that exceed the attachment of layer 2
      (still 8MM). I.e., it is now equivalent to a ``6MM in excess of 8MM''
      layer.
  }
  \label{fig:layered1}
  \vspace{5mm}
\end{figure}

\begin{figure}
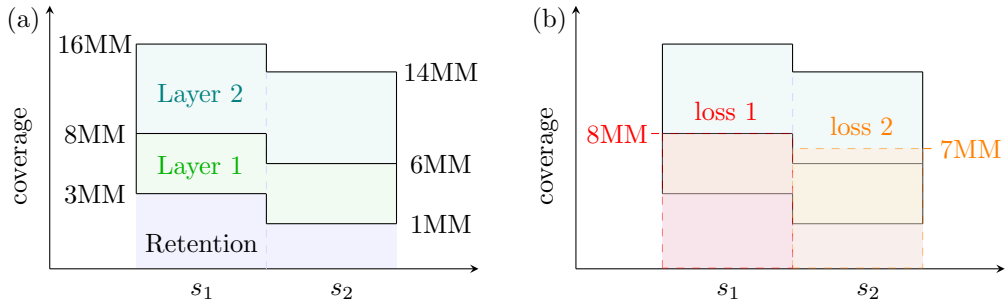

  \centering
  \include{figures/shifted}
  \vspace{-8mm}
  \captionsetup{width=.8\linewidth}
  \caption{
    \textbf{(a)} Multi-attachment reinsurance profile: Layer 1 has two peril
      subgroups ($s_{1}$ \& $s_{2}$) with respective attachments of 3MM \& 1MM
      and limit 5MM. Layer 2 has respective attachments of 8MM \& 6MM for
      $s_{1}$ \& $s_{2}$. Neither layer has a reinstatement clause.
    \textbf{(b)} Example: Distribution of 8MM loss 1 to perils $s_{1}$ (red)
      and 7MM loss 2 to perils $s_{2}$ (orange) between insurer and reinsurers:
      First, loss 1 occurs and the first 3MM of losses are retained by the
      insurer, then the next 5MM are ceded to the reinsurer of Layer 1.
      Secondly, loss 2 occurs and the first 1MM in losses are retained by the
      insurer. Layer 1 has been exhausted by loss 1 so cannot be used -- a
      further 5 MM which is not recovered. The remaining 1MM in losses between
      6MM and 7MM are ceded to the reinsurer of Layer 2.
  }
  \label{fig:layered2}
  \vspace{5mm}
\end{figure}

%%%%%%%%%%%%%%%%%%%%%%%%%%%%%%%%%%%%%%%%%%%%%%%%%%%%%%%%%%%%%%%%%%%%%%%%%%%%%%%
\vspace{-5mm}
\subsubsection*{Layer Structure:}
For each group $\group$, the reinsurance contract consists of one or multiple
layers. See Fig.~\ref{fig:layered1} for a schematic illustration of a simple
layer structure and Fig.~\ref{fig:layered2} for an attachment-shifted layer
structure. Each layer $\layer$ for a peril group $\group$ is characterized by:

\begin{itemize}
  \item The base attachment point: $\attach_{\layer} \geq 0$
  \item A subgroup-specific attachment shift: $\shift_{\subgroup}$
  \item The layer limit: $\limit_{\layer} > 0$
  \item The number of reinstatements: $\reinst_{\layer}$
\end{itemize}
The effective attachment point for a loss in subgroup $s$ is:
\begin{equation}
  \attach_{\layer,\subgroup} = \attach_{\layer} + \shift_{\subgroup}
  \label{eq:effective_attachment}
\end{equation}

%%%%%%%%%%%%%%%%%%%%%%%%%%%%%%%%%%%%%%%%%%%%%%%%%%%%%%%%%%%%%%%%%%%%%%%%%%%%%%%
\vspace{-5mm}
\subsubsection*{Recovery:} For an individual event with loss $\loss_\event$,
due a peril that is covered in subgroup $s$, the event loss recovery from layer
$\layer$ is:
\begin{equation}
  \recovery_{\layer,\event} = \min\left(
    \max\left[0, \loss_e - \attach_{\layer,\subgroup}\right],
    \limit_{\layer}
  \right)\,.
  \label{eq:event_recovery}
\end{equation}
In other words, an amount between $[0, \limit_{\layer}]$, depending on by how
much the loss $\loss_{\event}$ exceeds the relevant attachment point
$\attach_{\layer,\subgroup}$.

A layer's \emph{aggregate} recovery per trial year $\trial$ and layer $\layer$
is the sum of all its per-event recoveries, bounded by the layer's aggregate
limit $\limit_{\mathrm{agg}} = \limit (1+\reinst_{\layer})$ with up to
$\reinst_\layer$ reinstatements:
\begin{equation}
  \recovery_{\trial,\layer} = \min\left(
    (1+\reinst_\layer)\limit_\layer,
    \sum_{\event\in\trial} \recovery_{\layer,\event}
  \right)\,.
  \label{eq:yearly-recovery}
\end{equation}

The average aggregate recovery per layer can then be found:
\begin{equation}
  \bar{\recovery_{\layer}}
    = \frac{1}{|\{t\}|}\sum_{\trial} \recovery_{\trial,\layer}\,.
  \label{eq:yearly-recovery-layer}
\end{equation}

We calculate the retained loss in a given trial year $\trial$ as the
sum of the gross event losses minus the recovery from all layers:
\begin{equation}
  \mathcal{R}_{\trial} =
    \sum_{\event\in\trial} \loss_\event -
    \sum_\layer \recovery_{\trial,\layer}\,.
  \label{eq:yearly-loss}
\end{equation}
%

%%%%%%%%%%%%%%%%%%%%%%%%%%%%%%%%%%%%%%%%%%%%%%%%%%%%%%%%%%%%%%%%%%%%%%%%%%%%%%%
\vspace{-5mm}
\subsubsection*{Pricing Estimation:} The actual cost for placing a specific
reinsurance contract is challenging to estimate -- market sentiments can
fluctuate over short periods of time based on recent events. Pricing
information is also limited to few data points based on queries/feedback via a
broker for individual contract proposals. We model pricing based on the
expected value premium principle~\citep{Chi2011}, with the reinsurance premium
for layer $l$ given by:
\begin{equation}
  \premium^{\text{rein}}_\layer = (1+\rho) \bar{\recovery_{\layer}}\,,
  \label{eq:pricing-expected}
\end{equation}
which is the average aggregate recovery ($\bar{\recovery_{\layer}}$) with an
additional safety loading factor $\rho > 0$, which simplifies and emulates
reinsurance market conditions. When aiming to closely approximate current
market conditions, this market factor is substituted with peril-specific curves
$f(\cdot)$ describing the relation of the average aggregate recovery to the
reinsurance premium rate:
\begin{equation}
  \rol = \frac{\premium^{\text{rein}}_\layer}{\limit_\layer}
    := f_\layer(\lol)
    = f\left(\frac{\bar{\recovery_\layer}}{\limit_\layer}\right)\,,
  \label{eq:pricing-curve}
\end{equation}
Here, $\rol$ (rate on line) represents the reinsurance premium of the layer
divided by the limit $\limit_\layer$ of the layer, while $\lol$ (loss on line)
denotes the average aggregate recovery per layer divided by the same limit. The
curve function $f_\layer$ is expected to be strictly increasing ($f_\layer(x) >
x$) and have an initial value greater than zero ($f_\layer(0) =
\rol_{\layer,\mathrm{min}} > 0$). I.e., even a layer that has no expected loss
from a given data set will have a non-zero premium due to the minimum
rate-on-line $\rol_{\layer,\mathrm{min}}$. With this we obtain the
curve-estimated premium of:
\begin{equation}
  \premium^{\text{rein}}_\layer = L_\layer\cdot f
    \left(\frac{\bar{\recovery_\layer}}{\limit_\layer}\right)
  \label{eq:premium-curve}
\end{equation}
%

%%%%%%%%%%%%%%%%%%%%%%%%%%%%%%%%%%%%%%%%%%%%%%%%%%%%%%%%%%%%%%%%%%%%%%%%%%%%%%%
\vspace{-5mm}
\subsubsection*{Reinstatement Premium:} The concrete stipulations for
reinstatement pricing can differ for each contract. One pattern observed in the
real-world examples we used as the basis for this study are reinstatement costs
that are triggered automatically when the layer's limit needs to be replenished
-- even partially -- for a \emph{potential} subsequent event (which may or may
not occur). I.e., if an event $e$ has magnitude $A+L/2$, it recovers half the
layer's limit and also immediately triggers a reinstatement of this first half
that was used (assuming $\reinst_\layer>0$). After this, the next event has
again a per-event limit of $L$, but the \emph{aggregate} limit is $L/2$
consumed (and $50\%$ of one reinstatement has been paid so far). Assuming
that the premium for one full reinstatement is the same as the original premium
for the layer (referred to as ``$1@100$''), the total reinstatement cost can be
expressed as:
\begin{equation}
  \premium^{\text{rstm}}_{\trial,\layer} 
    = \premium^{\text{rein}}_\layer \cdot \frac{\min(
        \reinst_\layer\cdot\limit_\layer,\,\recovery_{\trial,\layer}
      )}{\limit_\layer}\,,
  \label{eq:premium-reinst}
\end{equation}
where the reinsurance premium is multiplied by a fraction which is a value
between $[0,\reinst_\layer]$, depending on the losses recovered in trial year
$\trial$. Note that the maximum aggregate loss considered for this calculation
is $\reinst\cdot\limit$ (as opposed to $(1+\reinst)\limit$ in
Eq.~\ref{eq:yearly-recovery}), because the \emph{last} $L$ ceded do not trigger
another reinstatement. While we focus on Eq.~\ref{eq:premium-reinst} for this
study, alternative terms for the reinstatement cost could be implemented by
decoupling $\premium_\layer$ and $\premium_{\reinst,\layer}$ and/or using an
alternative formula for either component.

Because the reinstatement cost for a layer is variable (i.e., it depends on the
losses ceded to the layer in a given trial year), it can be considered to be
part of the loss (or, equivalently, a reduction in the amount recovered). With
this, we obtain for the net loss
\begin{equation}
  \netloss_{\trial} =
    \sum_{\event\in\trial} \loss_\event -
    \sum_\layer\left[\recovery_{\trial,\layer} -
    \premium_{\trial,\layer}^{\text{rstm}}\right]\,.
  \label{eq:yearly-loss-reinst}
\end{equation}
which is used in the calculation of the average loss $\sum_t \frac{
  \netloss_{\trial}}{|\{t\}|}$, as well as the TVaR and AEP constraints
described below.

%%%%%%%%%%%%%%%%%%%%%%%%%%%%%%%%%%%%%%%%%%%%%%%%%%%%%%%%%%%%%%%%%%%%%%%%%%%%%%%
\vspace{-5mm}
\subsubsection*{Net Profit:}
Without reinsurance, the net profit $\pi_{\trial}$ in a given trial year
$\trial$ is obtained by subtracting the losses from events in that year from
the gross profit: $\pi_{\trial} = \Pi^{\text{gross}} - \sum_{\event\in\trial}
\loss_e$. The gross profit ($\Pi^{\text{gross}}$) is a pre-computed constant
based on insurance premiums and costs, and for simplicity, we assume it remains
unchanged across different trial years. When reinsurance is applied, some of
these losses are recovered (Eq.\ref{eq:yearly-loss}), but the reinsurance
premium (Eq.~\ref{eq:premium-curve}) and reinstatement costs
(Eq.~\ref{eq:premium-reinst}) need to be paid:
\begin{equation}
  \pi_{\trial}
    = \grossprofit - \loss_{\trial} - \sum_\layer\left[
        \premium^{\text{rein}}_{\layer} + \recovery_{\trial,\layer}
        - \premium_{\trial,\layer}^{\text{rstm}}
      \right]
    \equiv \grossprofit - \netloss_{\trial}
      - \sum_{\layer}\premium^{\text{rein}}_\layer
  \label{eq:netprofit}
\end{equation}
%

%%%%%%%%%%%%%%%%%%%%%%%%%%%%%%%%%%%%%%%%%%%%%%%%%%%%%%%%%%%%%%%%%%%%%%%%%%%%%%%

\subsubsection*{Risk:}
\label{sec:tvar}
The value of a reinsurance contract comes from reducing an insurers risk
exposure. Here we use the Tail Value at Risk (TVaR) as our risk measure as it
focuses on the size and severity of extreme losses. It is the average over the
yearly net loss values which are above a percentile ($\mathbf{P_{\beta}}$):
\begin{equation}
  \mathrm{TVaR}_{\beta}[\netloss_{\trial}] = |\netloss_{\trial}| :
    \netloss_{\trial} > \mathbf{P}_{\beta}(
      \netloss_{\trial})\,.
  \label{eq:tvar}
\end{equation}

For TVaR, $\beta$ is typically 200 year percentile or $99.5\%$.

%%%%%%%%%%%%%%%%%%%%%%%%%%%%%%%%%%%%%%%%%%%%%%%%%%%%%%%%%%%%%%%%%%%%%%%%%%%%%%%
\vspace{-5mm}
\subsubsection*{Constraints:}
\label{sec:constraints}
The selection of potential reinsurance contracts is subject to constraints,
either due to regulatory requirements or business-imposed. We can transform
these constraints into mathematical inequality constraints and we consider the
following properties, along with threshold values for each. These thresholds
are called alert values and are set by the business:
\begin{itemize}

  \item {\it Aggregate Exceedance Probability:} 
    AEP measures the size of the \textit{total} accumulated net losses from all
    events in a given trial year. It considers the combined impact of multiple
    events, making it particularly valuable for assessing the annual portfolio
    risk. It is typically defined as a percentile of the net losses:
    \begin{equation}
      \mathrm{AEP}_{\beta} =
      \mathbf{P}_{\beta}(\netloss_{\trial})
      <\mathrm{AEP}_{\beta}^{\mathrm{alert}}\,.
      \label{eq:aep}
    \end{equation}
    The AEP is calculated over all perils and typically looks at 5, 10, and
    50 year percentile $\beta$.

  \vspace{2mm}
  \item {\it Occurrence Exceedance Probability}:
    OEP measures the probability that at least one single event in a given time
    period will exceed a specified loss threshold. It focuses on
    \textit{individual} event severity rather than cumulative impact. OEP is
    particularly useful for analyzing catastrophic event exposure and setting
    single-event coverage limits. It is defined as a percentile of a
    distribution made up of the single largest event (net of reinsurance) of a
    given peril $p$ within a year. To this end, we select the largest event in
    each year and apply the relevant reinsurance layers to its loss:
    $\loss_{\peril,\trial}^{\mathrm{max}} = \max_{\event\in\peril,\trial}
      (\loss_\event-\sum_\layer\recovery_{\layer,\event})$:
    \begin{equation}
      \mathrm{OEP}_{\beta, \peril} =
      \mathbf{P}_{\beta}(\loss_{\peril,\trial}^{\mathrm{max}})
      <\mathrm{OEP}_{\beta,\peril}^{\mathrm{alert}}\,.
      \label{eq:oep}
    \end{equation}
    The OEP is calculated over specific perils, usually peak perils in North
    America and Europe and generally its $\beta$ is 200 year.

  \vspace{2mm}
  \item {\it Probability of attaching:} The probability of ceding a non-zero
    amount to a given reinsurance layer is given by:
    \begin{equation}
        P_{\layer}^{\mathrm{att}}
        = \frac{|\{\recovery_{\trial,\layer} > 0\}|}{|\{t\}|}
        < P_{\layer}^{\mathrm{att},\mathrm{alert}}\,.
      \label{eq:attach_prob}
    \end{equation}
    This particular metric is of interest because a high probability of
    attaching may be a hindrance to placement of the reinsurance contract.

\end{itemize}

%%%%%%%%%%%%%%%%%%%%%%%%%%%%%%%%%%%%%%%%%%%%%%%%%%%%%%%%%%%%%%%%%%%%%%%%%%%%%%%
%%%%%%%%%%%%%%%%%%%%%%%%%%%%%%%%%%%%%%%%%%%%%%%%%%%%%%%%%%%%%%%%%%%%%%%%%%%%%%%
\section{Markov-Chain Optimization}
\label{sec:markov-chain}

We analyze the prospects of simulated annealing as a search method to
efficiently identify high-quality reinsurance structures. Given sufficient
time, frequent restarts, and appropriate parametrization, this approach has a
high probability of identifying the global optimum even in a large search
space. In particular, simulated annealing allows us to identity near-optimal
solutions in a much shorter time frame and has the additional ability to “warm
start” from a predefined set of initial solutions to focus on a specific region
of the solution space.
% Optimization of an objective function using such a Markov chain process
% requires efficient evaluation of individual contracts and a suitable
% neighborhood structure that allows traversal to similar contract structures
% at each step.
In what follows, we outline a Markov-chain-based framework that integrates
efficient evaluation of individual contracts with a custom neighborhood
structure that allows for traversal to similar contract structures at each
step. 

%%%%%%%%%%%%%%%%%%%%%%%%%%%%%%%%%%%%%%%%%%%%%%%%%%%%%%%%%%%%%%%%%%%%%%%%%%%%%%%
\vspace{-5mm}
\subsection{Preprocessing} We transform simulated loss data into a format that
allows for rapid computation of recoveries for layers. For each peril $p$ and
trial year $\trial$, we aggregate simulated losses into a function
$D_{\trial,\peril}(x)$ defined as the sum of all losses below a given threshold
$x$:
\begin{equation}
    D_{\trial,\peril}(x) = \sum_{\event \in \trial,\peril} \min(x, \loss_e)\,.
\end{equation}
This transformation allows for the calculation of the (unconstrained) recovery
as the difference:
\begin{equation}
  \recovery_{\trial,\layer,\peril}' =
    D_{\trial,\peril}(\attach_{\layer} + \limit_{\layer}) -
    D_{\trial,\peril}(\attach_{\layer})\,.
  \label{eq:cumsum-diff}
\end{equation}
For a limited set of possible attachment points $\{x\}$ (i.e., artificially
discretized or constrained to round numbers), the values of $D_{\trial,\peril}$
at these points can be pre-calculated and stored in a database. If the peril
subgroups of interest are known a priori, this transformation can also be done
at the subgroup level. The allows us to leverage grouped-sum database queries
to efficiently recall the recovered amount $\recovery$ for combinations of
perils and multiple trial years, which is crucial for the optimization
algorithms described in subsequent sections.

%%%%%%%%%%%%%%%%%%%%%%%%%%%%%%%%%%%%%%%%%%%%%%%%%%%%%%%%%%%%%%%%%%%%%%%%%%%%%%%
\vspace{-5mm}
\subsection{Simulated Annealing}
\label{sec:simulated-annealing}

We employ a Markov chain optimization approach based on simulated annealing.
This method constructs a chain of contracts that explores the parameter space
of possible reinsurance structures via small variations at each step. This
transition from one contract structure $C$ to the next is guided by improving
the objective function and occasionally accepting suboptimal moves to escape
local optima. Formally, a Markov chain is defined by:
\begin{equation}
  P(C_{n+1} = x \,|\, C_1 \dots, C_n) = P(C_{n+1} \,|\, C_n)\,.
\end{equation}
That is, the probability of picking each neighboring state as the next depends
only on the current one (not the full chain). Here we consider the state space
of reinsurance contracts, defined by:
\begin{equation}
  C = \{(\attach_{\layer,\subgroup}, \limit_{\layer},
         \reinst_{\layer}, \group, \subgroup, \Delta_{\subgroup}),\ldots\}\,,
\end{equation}
which is a collection of reinsurance layers (each defined by their attachment
$\attach_{\layer,\subgroup}$, limit $\limit_{\layer}$, reinstatements
$\reinst_{\layer}$) and the perils to which they apply (defined by the peril
groups $\group$, subgroups $\subgroup$ and subgroup shifts
$\Delta_{\subgroup}$). We impose constraints to the state space to ensure that
states remain within realistic boundaries for a reinsurance contract. These
constraints vary according to the specific business requirements of each
contract.

%%%%%%%%%%%%%%%%%%%%%%%%%%%%%%%%%%%%%%%%%%%%%%%%%%%%%%%%%%%%%%%%%%%%%%%%%%%%%%%
\vspace{-5mm}
\subsubsection*{Acceptance Probability:}
Transitions in contract space are accepted according to the Metropolis
criterion \citep{Metropolis1953}:
\begin{equation}
  P_{\mathrm{accept}} = \min\left\{
    1, \exp\left[(\objective(C') - \objective(C))/T\right]
  \right\}
\end{equation}
where $\objective(C)$ is the objective function that we wish to maximize, $C
\to C'$ are the contract structures transitioned from (to) and $T>0$ is the
simulated annealing temperature. This temperature parameter $T$ controls the
probability of accepting transitions that worsen the the objective function;
potentially allowing it to escape from local maxima in $\objective(C)$. Moves
that improve the objective function are always accepted ($\objective(C') >
\objective(C)\Rightarrow P_{\mathrm{accept}} = 1$).

The effectiveness of simulated annealing depends heavily on this temperature
parameter: As it decreases from an initial high temperature to a lower one, the
algorithm transitions from exploration (accepting worse positions to move
around the solution space) to exploitation (identifying the best solution in a
local neighborhood). The geometric cooling scheme is one of the most widely
adopted approaches \citep{VanLaarhoven1992}. Under this scheme, the temperature
$T$ at step $k$ is given by:
\begin{equation}
  T_k = T_0 \alpha^k := T_0 \left(
    \frac{T_f}{T_0}
  \right)^\frac{k}{k_\mathrm{max}}
\end{equation}
where $T_0$ is the initial temperature and $\alpha$ is the cooling rate
(typically $0.8 < \alpha < 0.99$). The initial temperature $T_0$ should be
sufficiently high to allow for a wide exploration of the solution space, while
the final temperature $T_f = T_0 \alpha^{k_\mathrm{max}}$ should drive
exploitation: Focusing on local improvements while only rarely accepting
worsening transitions. We calculate $\alpha$ based on $T_0$, $T_f$ and
$k_{\mathrm{max}}$, which we use as input parameters.

%%%%%%%%%%%%%%%%%%%%%%%%%%%%%%%%%%%%%%%%%%%%%%%%%%%%%%%%%%%%%%%%%%%%%%%%%%%%%%%
\vspace{-5mm}
\subsection{Implementation \& Results}
Our implementation relies on pre-computed cumulative loss sums stored in an SQL
database to efficiently calculate expected losses. By storing cumulative sum
functions $D_{\trial,\peril}(x) = \sum_{\event\in\trial}\min(x,\loss_\event)$,
we can compute the losses of a peril that are applicable to a layer through the
difference operation $D_{\trial,\peril}(\attach+\limit) -
D_{\trial,\peril}(\attach)$ (see Eq.~\ref{eq:cumsum-diff}). This approach
significantly reduces the computational overhead for calculating the layer
recovery with aggregate limit
\begin{equation}
  \recovery_{\trial,\layer} =
    \min\left[
      (1+\reinst)L,\, \sum_{\peril\in\layer}\left(
        D_{\trial,\peril}(\attach_\layer+\limit_\layer)-
        D_{\trial,\peril}(\attach_\layer)\right)\right]\,.
  \label{eq:recovery-via-diff}
\end{equation}
Storing the full information that allows evaluating $D_{\trial,\peril}(x)$ at
any value of $x$ requires storage comparable to the original dataset (the
individual event losses could be recovered through sufficient queries).
However, in practice, reinsurance contracts are drafted with round numbers for
the attachments $\attach_\layer$ and limits $\limit_\layer$. This effectively
discretizes the allowed values and allows us to pre-compute and store the value
of $D_{\trial,\peril}(x)$ only at this fixed set of values for $x$. For our
dataset of $50,000$ trial years with $42,141,370$ events, discretization to 80
round values reduces the size from 2.4GB (hdf5, individual event losses), to
33MB (SQL database, cumulative loss values) -- a reduction of 74$\times$. If
the grouping is known a priori (or constrained to a specific set of allowed
groupings), we can instead store $D_{\group,t}$ to avoid calculating the sum of
perils in Eq.~\ref{eq:recovery-via-diff}.

%%%%%%%%%%%%%%%%%%%%%%%%%%%%%%%%%%%%%%%%%%%%%%%%%%%%%%%%%%%%%%%%%%%%%%%%%%%%%%%
\vspace{-5mm}
\subsubsection*{Decision Variables}

The decision variables defining a contract $C$ can be inferred per contract
layer structure as detailed in Sec.~\ref{sec:model}.

%%%%%%%%%%%%%%%%%%%%%%%%%%%%%%%%%%%%%%%%%%%%%%%%%%%%%%%%%%%%%%%%%%%%%%%%%%%%%%%
\vspace{-5mm}
\subsubsection*{Neighborhood Structure}

\begin{figure}
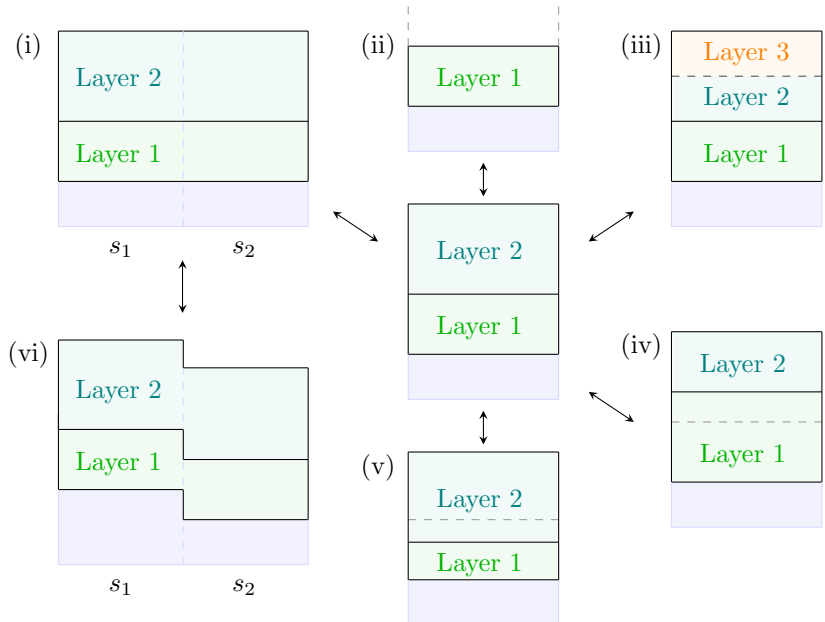

  \vspace{3mm}
  \centering
  \include{figures/moves}
  \vspace{-1mm}
  \captionsetup{width=.8\linewidth}
  \caption{
    Visualization of the different Markov chain moves from an initial state
    (center: one tower with two layers), and their reversal (to the center).
    \textbf{(i)} Creation and removal of a subgroup.
    \textbf{(ii)} Removal and addition of a layer.
    \textbf{(iii)} Splitting and joining of a layer.
    \textbf{(iv)} Adjustment of layer boundary: Layer 1 grows while Layer 2
      shrinks.
    \textbf{(v)} Adjustment of layer boundary with shift above: Layer 2 retains
      the same limit.
    \textbf{(vi)} Subgroup shift from (i).
    %Not depicted: Changing of reinstatement terms.
  }
  \label{fig:moves}
  %\vspace{-5mm}
\end{figure}

The neighborhood of each contract $C$ is defined by permissible variable
adjustments. As schematically shown in Fig.~\ref{fig:moves}, in our model moves
between states are generated by perturbing individual variables as follows:
\begin{enumerate}
  \item Adjusting groups and subgroups.
  \vspace{1mm}
  \item Adding/removing a layer.
  \vspace{1mm}
  \item Joining or splitting existing layers.
  \vspace{1mm}
  \item Adjusting attachment points
    ${\attach}_{\layer,\subgroup}$ or limits ${\limit}_{\layer}$.
  \vspace{1mm}
  \item Adjusting a limit while shifting everything above in tower by the same
    amount.
  \vspace{1mm}
  \item Changing subgroup shifts ${\Delta}_{\subgroup}$.
  \vspace{1mm}
  \item Changing reinstatement terms ${\reinst}_{\layer}$.
\end{enumerate}

In accordance with the sample contracts studied in this work, we assume that
the reinsurance layers form ``towers'': The upper boundary ($\attach + \limit$)
of one layer coincides with the attachment of the layer above it. That is,
there are no gaps in the coverage above the retention. Peril grouping and
subgrouping is assumed to be identical for each layer in the tower, which
simplifies the adjustment of groups~(i).

Adding or removing a layer~(ii) provides the primary means of changing the
structure of a contract beyond adjusting attachments and limits. Only the top
and bottom layers in a tower can be chosen for this move. These moves tend to
cause a large change in the objective function and are typically only accepted
during initial solution building or when they directly counter-act a current
constraint violation.

Joining, splitting (iii) and changing attachment/limit (iv) serves the purpose
of fine-tuning a solution during the final exploitation phase of the
optimization, while adjustment with shift (v) allows a layer to grow/shrink
without changing the size of the layers above it. This is particularly useful
when a minimum layer size is enforced: Small layers would prevent growing the
layer below it otherwise. We are not considering the symmetric move with
shifting the layers \emph{below}, because changes to the lower layers have a
larger impact on the objective function.

Adjustment of subgroup shifts~(vi) allows for a different retention for perils
with different levels of activity which are in the same group (but different
subgroups). Finally, changing the reinstatement terms (vii) typically has a
small effect on the objective function, but can help in alleviating violated
constraints.

%%%%%%%%%%%%%%%%%%%%%%%%%%%%%%%%%%%%%%%%%%%%%%%%%%%%%%%%%%%%%%%%%%%%%%%%%%%%%%%
\vspace{-5mm}
\subsubsection*{Objective Function:} We use the expected net profit
(Eq.~\ref{eq:netprofit}) as the primary objective function component in the
optimization (with additional penalty terms $\lambda_K$ for constraints).
\begin{equation}
  \objective(C) = \mathrm{avg}_\trial(\netprofit{\trial}(C)) + \sum_K \lambda_K(C)\,,
  \label{eq:objective-main}
\end{equation}
Its evaluation requires the expected loss for each layer -- both to compute the
expected recovered amount (Eq.~\ref{eq:yearly-recovery}) (which reduces the
loss) and the premium and reinstatement costs (Eqs.~\ref{eq:premium-curve}
and~\ref{eq:premium-reinst}). Identifying the yearly losses via
Eq.~\ref{eq:recovery-via-diff} can be formulated as a grouped SQL query, which
allows us to leverage a database engine that is optimized for these kinds of
lookups. Additionally, we employ multi-tiered LRU-caching at the level of the
individual layer and for full contracts to limit the effective number of
database queries. This proves particularly effective since the Markov chain
explores similar contracts in sequence. I.e., most transitions in our defined
neighborhood structure only affect an individual layer or stack of layers --
calculation results for all unmodified parts of the contract can be found in
the cache.

%%%%%%%%%%%%%%%%%%%%%%%%%%%%%%%%%%%%%%%%%%%%%%%%%%%%%%%%%%%%%%%%%%%%%%%%%%%%%%%
\vspace{-5mm}
\subsubsection*{Constraints:} We consider the constraints described in
section~\ref{sec:model} along with a threshold value they should not exceed.
These are optional and implemented as soft constraints: For every constraint
$K$, such as TVaR, with threshold value $\tau_K$ set, we add a penalty term of
the form
\begin{equation}
  \label{eq:penalty-form}
  \lambda_K = -\gamma_K \cdot \mathrm{ReLU}(K-\tau_K)\,,
\end{equation}
where $\mathrm{ReLU}(x)=\max(0,x)$ is the rectified linear unit function and
$\gamma_K$ is a scaling prefactor. This has the effect of reducing the
objective function linearly when a constraint is violated, but still allowing
exploration of the solution space via ``small'' violations. In particular, this
also solves the starting position problem: It allows us to start from an
initially empty contract (i.e., no reinsurance layers) -- which almost
certainly violates some constraints; otherwise it would be the highest
net-profit solution. To ensure that these soft violations do not surface as
solutions, additional filtering is applied to the values stored as best
candidates during the optimization to remove any that violate any of the
constraints.

Most constraints can be evaluated from the data used to compute the net profit.
However, for the OEP constraint (Eq.~\ref{eq:oep}), information about the
individual largest event for each peril and trial year is not available from
$D_{\trial,\peril}(x)$. In order to still be able to evaluate it, we separately
store the largest-in-year events for each peril.

%%%%%%%%%%%%%%%%%%%%%%%%%%%%%%%%%%%%%%%%%%%%%%%%%%%%%%%%%%%%%%%%%%%%%%%%%%%%%%%
\vspace{-5mm}
\subsubsection*{Optimization:} We use simulated annealing with restarts. Each
run can be started either from an initially empty contract or a (set of)
initial contract(s). When starting from an empty contract, we increase the
probability of proposing a layer addition during the initial exploration phase.
This bias enables initial solution building and speeds up the time to reach
parts of the feasible solution space. After this, regular exploration and
eventually exploitation with simulated annealing resumes according to the
temperature schedule. Throughout the process, we keep track of the most
promising solutions (as evaluated by the objective function) and return it
either at the end or interactively throughout the process. The latter enables
visual feedback on the progress of the optimization and contracts identified so
far.

\begin{figure}
  \begin{center}
  \includegraphics[width=0.6\textwidth]{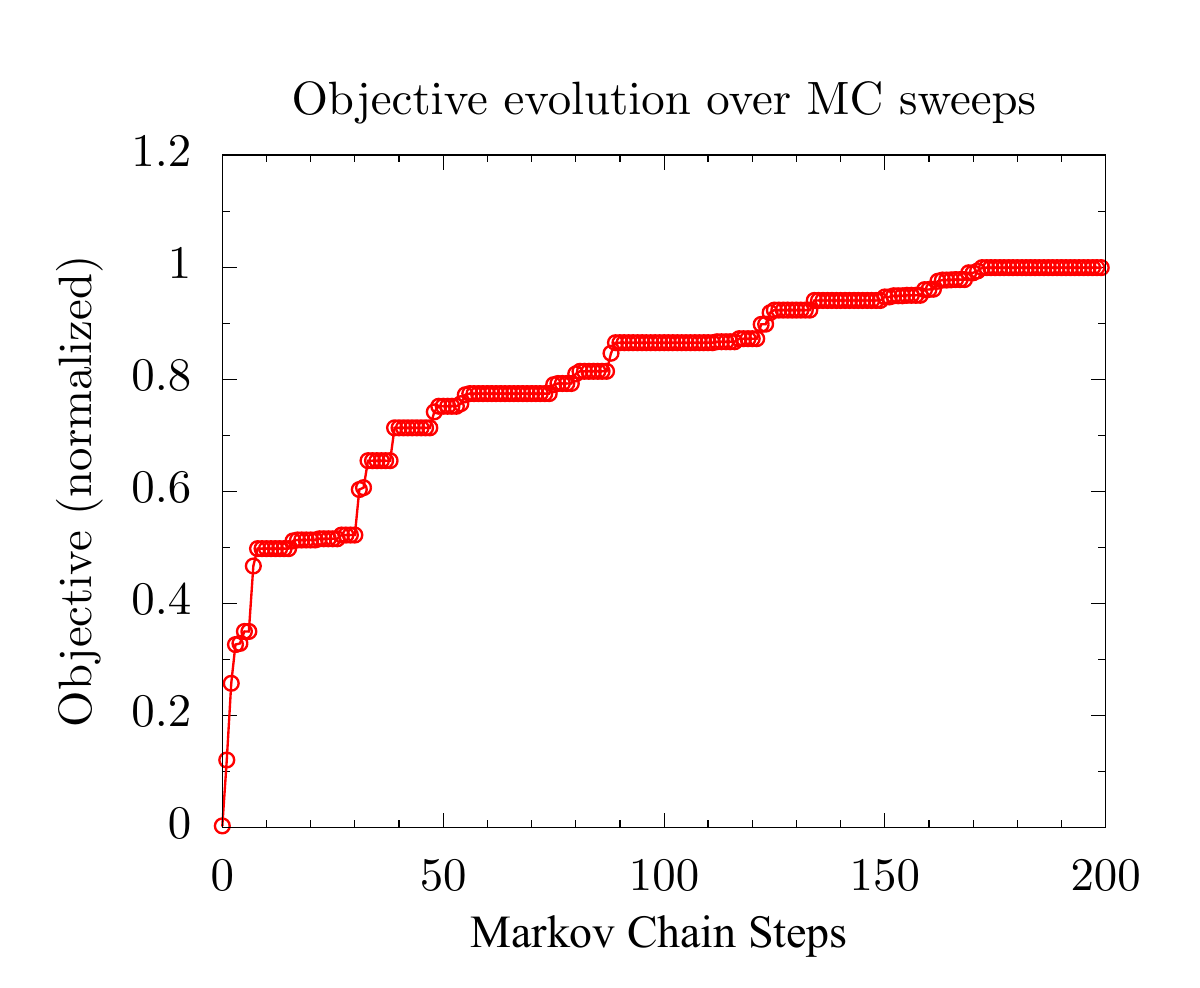}
  %\gr
  \captionsetup{width=0.8\textwidth}
  \caption{
    While substantial improvements in the objective function can typically be
    found early on (depending on the starting point), these increments become
    smaller as the optimizer approaches the best solution found. While
    simulated annealing is a heuristic and cannot provide a guarantee about
    optimality, we can increase the probability of finding a solution that is
    close to optimal by (i) ensuring sampling in each Markov chain is long
    enough to reach this region of small increments and (ii) running several
    Markov chains in parallel (with different random seeds).
  }
  \label{fig:improvement}
  \end{center}
\end{figure}

\begin{figure}
  \begin{center}
  \includegraphics{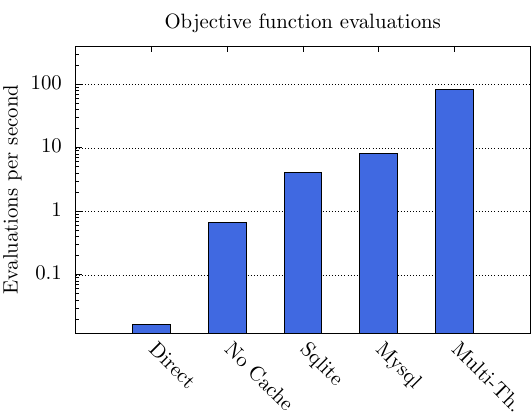}
  \captionsetup{width=0.8\textwidth}
  \caption{
    Rapid objective function evaluation is key to a numerical optimization
    approach. When computed directly from the full event data set, this takes
    several seconds for each evaluation. Dataset normalization and cumulative
    sum representation reduce this time to below one second. With caching, a
    capable database backend and multi-threading, we can evaluate more than 80
    contracts per second (Ryzen 9 5950X, 96GB ram).
  }
  \label{fig:evaluation}
  \end{center}
\end{figure}

A visualization of this improvement over time (represented by an increase in
the objective function) is shown in Fig.~\ref{fig:improvement}. The sample
optimization run is for a contract encompassing $43$ perils with partially
restricted grouping on the full dataset of $50,\!000$ trial years. The
optimization considered $32,\!000$ alteration steps in the Markov chain within
$382$ seconds. The main contributor to the runtime of the optimization is the
evaluation of the objective function -- in particular calculating the amount
recovered for each layer and trial year. Fig.~\ref{fig:evaluation} indicates
the number of objective function evaluations achieved per second for different
architectural setups, highlighting the importance of leveraging a performant
database and caching to improve the sampling speed.

\begin{figure}
  \begin{center}
  \includegraphics{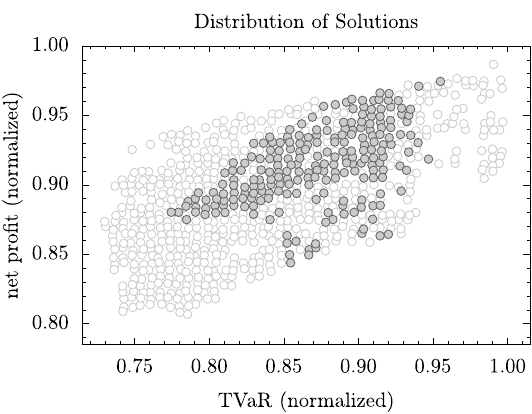}
  \captionsetup{width=0.8\textwidth}
  \caption{
    Distribution of contracts visited as part of an extended optimization run
    with a focus on exploration. The dots were randomly sampled from all
    solutions considered. Dark gray dots indicate a (feasible) contract that
    satisfies the constraints, while a light one does not. The efficient
    frontier in the net-profit vs TVaR trade-off is the top left boundary of
    the area with gray dots, where profit is maximized for a given risk.
    (Both axes are normalized w.r.t. their value without reinsurance.)
  }
  \label{fig:space}
  \end{center}
\end{figure}

This rapid exploration of potential contracts allows us to render the
distribution of solutions in, e.g., the net profit vs TVaR dimensions
(Fig.~\ref{fig:space}). Along with filtering for solutions that meet all
constraints, we can visualize the efficient frontier of contracts (upper-left
boundary). Empirically we find that -- unless dictated otherwise by specific
constraints -- a single large layer is better than multiple smaller ones and
grouping many perils is beneficial; especially since subgroup shifts allow
adjusting the retention for perils in subgroups.

%%%%%%%%%%%%%%%%%%%%%%%%%%%%%%%%%%%%%%%%%%%%%%%%%%%%%%%%%%%%%%%%%%%%%%%%%%%%%%%%
%%%%%%%%%%%%%%%%%%%%%%%%%%%%%%%%%%%%%%%%%%%%%%%%%%%%%%%%%%%%%%%%%%%%%%%%%%%%%%%%
\section{Quantum Branch \& Bound}
\label{sec:branch-bound}

In this work, we also investigate the use of branch and bound algorithms to
optimize the structure of reinsurance contracts. These are exact algorithms
that guarantee an optimal solution. Specifically, we focus on quantum branch \&
bound (QBB) which has the potential to leverage future fault-tolerant quantum
computers.

Branch \& bound (B\&B) is a widely used method for solving constrained
optimization problems. It focuses on finding the minimum-cost solution which
satisfies all constraints. Two key components are involved: a bounding function
that returns a cost lower bound of a subset of possible solutions and a
branching rule that divides a subset of possible solutions into smaller
subsets. The algorithm eliminates subsets of solutions where the cost lower
bound is higher than the current best cost, eliminating the need to explore
every possible option. If the lower bound of a subset of solutions is less than
or equal to the current best cost, it is divided into smaller, disjoint subsets
using the branching rule. This approach also applies to maximization problems.
It is evident that the key in the algorithm's performance is a tight bounding
function with an associated branching rule as otherwise the algorithm would
exhaustively evaluate each possible option.

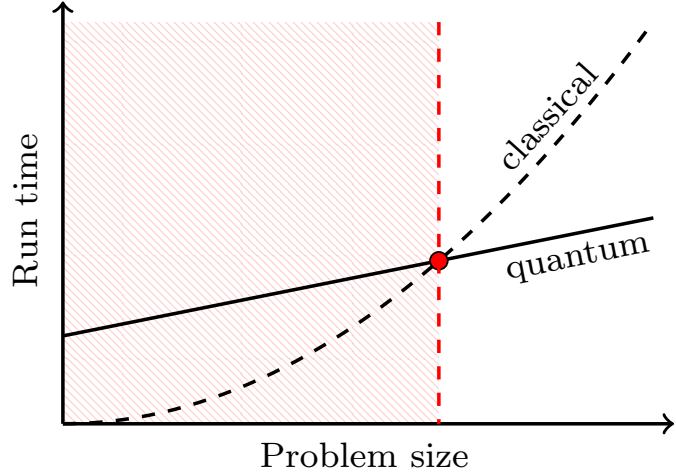
\begin{figure}[t]
  \begin{center}
    \resizebox{.6\linewidth}{!}{
      \begin{tikzpicture}
        \fill[pattern=north west lines,pattern color=red!20] (0,0) rectangle (2.99253,3.2);
        
        \node (x1) at (5,0) {};
        \node (y1) at (0,3.5) {};
        
        \draw[->,thick] (0,0) -- (x1);
        \draw[->,thick] (0,0) -- (y1);
        
        \node (ylabel) at (-0.3,1.7) {\rotatebox{90}{\footnotesize Run time}};
        \node (ylabel) at (2.4,-.25) {\rotatebox{0}{\footnotesize Problem size}};
        
        \draw[thick,dashed,domain=0:4.7,smooth,variable=\x] plot ({\x},{0.145*\x*\x});
        \draw[thick,domain=0:4.7,smooth,variable=\x] plot ({\x},{0.2*\x+0.7});
        
        \node at (4.1,1.32) {\rotatebox{11.2}{\footnotesize quantum}};
        \node at (3.85,2.45) {\rotatebox{48}{\footnotesize classical}};
        
        \draw[draw,fill=red] (2.99253,0.2*2.99253+0.7) circle (2pt);
        \draw[dashed,red,thick]  (2.99253,3.2) --  (2.99253,0);
      \end{tikzpicture}
    }
    \captionsetup{width=.8\linewidth}
    \caption{
      Illustration of a typically scenario with a quantum advantage.
      While the quantum algorithm scales better than the classical
      counterpart, the high overheads of quantum computers (quantum error
      correction and fault tolerance) negate a quantum speedup for small
      problem sizes.
      It is important to determine the \textit{crossover point},  i.e.,
      the problem size after which the quantum algorithm outperforms its
      classical counterpart. Figure copied from
      \citep{steiger2019advantages}.
    }
    \label{fig:crossoverpoint}
  \end{center}
\end{figure}

QBB achieves a nearly quadratic speedup over classical B\&B in terms of the
number of calls to the branching and bounding operators
\citep{montanaro2020quantum}. A scaling advantage for a quantum algorithm is a
necessary but not a sufficient condition for a runtime advantage of a quantum
computer for solving problems in practice. Given a scaling advantage, it is
clear that there exists a problem size at which the quantum algorithm starts
outperforming the classical algorithm, see Fig.~\ref{fig:crossoverpoint} for an
illustration. However, this quantum advantage is only relevant in practice if
the runtime at that cross-over point is reasonable.
For our analysis, we allow a maximum runtime of $10^6$ seconds (approximately
11 and a half days) for practical reasons such as potential changes in the
pricing of contracts during the runtime of the algorithm and to enable
reasonable turnaround times during contract negotiations between reinsurers and
their clients.

In this setting, it is crucial to examine not just the asymptotic scaling of
the algorithm, but also the overheads associated with fault-tolerant quantum
computing
\citep{steiger2019advantages,Babbush2021QuadraticSpeedup,hoefler2023disentangling}.
To this end, we adopt a simple model for a quadratic quantum
speedup~\citep{Babbush2021QuadraticSpeedup}, where the runtime of the classical
algorithm is given by
\begin{equation}\label{eq:t_classical}
  T_{c}(N) = t_c \cdot M \cdot N
\end{equation}
and the runtime of the quantum algorithm is given by
\begin{equation}\label{eq:t_quantum}
  T_{q}(N) = t_q \cdot M \cdot \sqrt{N},
\end{equation}
where $N$ is the number of calls to the classical oracle, $M$ is the number of
elementary operations per oracle call (assumed to be the same for the quantum
and classical implementation), and the times $t_{c}$, and $t_{q}$ are the time
per operation on a classical and on a quantum computer, respectively. This
simplification is favorable toward quantum computers, as implementing classical
oracles would introduce additional overheads, e.g., due to reversibility. 

In the case of QBB, $N$ is the number of calls to the branching and bounding
functions, and $M$ is the number of elementary gate operations (e.g., AND
gates) required to implement these two functions. To achieve a practical
quantum advantage, we must have $T_q < T_c$, from which it follows that the
B\&B tree size $N$ must be larger than $\left(t_q/t_c\right)^2$ and
additionally $T_q < 10^6 s$ (11.5 days), which limits the number of elementary
operations $M$ allowed per oracle call~\citep{hoefler2023disentangling}.

As an example of the $t_q/t_c$ ratio, let us consider the classical AND
operation and compare this to the quantum Toffoli gate, which is the quantum
analog of the classical AND gate~\citep{barenco1995elementary}. Assuming
superconducting qubits and the surface code for achieving fault tolerance
\citep{kitaev2003fault, Fowler2012SurfaceCode, gidney2019efficient,
Babbush2021QuadraticSpeedup} estimate the time for a Toffoli gate to be
$\approx 170 \mu s$. In contrast, a classical CPU running at 3GHz can implement
multiple AND gates in a single clock cycle ($\frac 13 ns$), from which we see
that $t_q/t_c> 5.1\cdot 10^5$ for AND gates. Furthermore, for the comparison to
be fair, one would also have to take into account the potential parallelism in
both cases, as well as the classical hardware required for quantum error
correction.

In this work, we investigate whether there is a medium-term quantum advantage
of QBB for solving a simplified version of the reinsurance optimization problem
previously described. Given that the ratio $t_q/t_c$ is large, it is crucial to
limit the number of operations required for the branching and cost function. We
thus focus on finding an inexpensive, but effective application-specific bound
for the simplified optimization problem.

\subsection{Simplified reinsurance problem}

For the analysis of a branch \& bound approach to identifying the optimal
layers, we consider a simplified problem: We assume that we are given a set of
$n$ peril groups, each peril group consists of one or more perils $p$, and each
peril group has one excess-of-loss layer $\layer_g$. Where $g=1,...,n$ and the
layers are specified by the variables $\layer_1,\dots,\layer_n$ each comprised
of an attachment point $\attach_g$ and a limit $\limit_g$. The reinsurance
premium and the net profit are estimated using the expected value premium
principle (according to Eqs.~\ref{eq:pricing-expected} and~\ref{eq:netprofit}).
This simplifies the optimization problem as we only need to find an optimal
attachment point and limit for each peril. In a two-stage optimization process,
one could then divide this large single layer into multiple smaller ones to
optimize further.

%%%%%%%%%%%%%%%%%%%%%%%%%%%%%%%%%%%%%%%%%%%%%%%%%%%%%%%%%%%%%%%%%%%%%%%%%%%%%%%%
\vspace{-5mm}
\subsubsection*{Risk constraints:}
Given the market factor $\rho$ in Eq.~\ref{eq:pricing-expected}, the
reinsurance premium is higher than the recovery from a layer. However,
reinsurance is necessary for insurance companies to reduce the risks of large
losses in any given year based on various risk measures. We start with a
general risk constraint $K$ of the following form:
\begin{equation}\label{eq:risk-constraint}
  K(\layer_1, \dots, \layer_n) \leq K_{\max},
\end{equation}
where $\layer_1, \dots, \layer_n$ are the single layers of each peril group and
the threshold $K_{\max}$ is a constant. Examples of $K(\cdot)$ might be the
tail value at risk (Eq.~\ref{eq:tvar}), aggregate exceedance probability
(Eq.~\ref{eq:aep}) and occurrence exceedance probability (Eq.~\ref{eq:oep}).

\subsection{Synthetic data generation and compression}
\label{sec:data-compression-generation}

Loading data into quantum computers is known to be very
expensive~\citep{hoefler2023disentangling,jaques2023qram}. We therefore aim to
reduce the number of data points to be loaded. To compress the data, we use the
attachment probability constraint (Eq.~\ref{eq:attach_prob}). Enforcing a
strict maximum attachment probability of $0.1$ (every 10 years) allows us to
compute a minimum attachment for each layer $\attach^{\min}_{g}$. Event losses
below this minimum attachment can be summed up as a total loss for each year
that is never ceded to a reinsurer. Only events with losses greater than this
minimum attachment need to be loaded. 
Specifically, for our compression experiments, we use a synthetic data set
derived from a real-world data set, consisting of 34 perils with events and a
total of 42,141,370 events over 50,000 simulated years. Enforcing the maximum
attachment probability constraint allows the number of events to be reduced by
a factor of $107\times$ to 394,067 events.

\newcommand{\pg}{\peril}
\newcommand{\AL}{\attach^{\mathrm{min}}_{g}}
\newcommand{\AU}{\attach^{\mathrm{max}}_{g}}
\newcommand{\LU}{\limit^{\mathrm{max}}_{g}}

We have also implemented this compression algorithm as part of the classical
B\&B code (discussed below in Section~\ref{sec:classical-bb}). The maximum
attachment probability $P_{\mathrm{attach}}\in[0,1]$ is a user-supplied command
line argument that is used to (i) compute minimum attachment values $\AL$ for
each peril group $g$ and (ii) compress the generated synthetic data. For each
peril group $g$, the B\&B algorithm then has to find the optimal attachment
$\attach_g\in [\AL, \AU]$ and limit $\limit_g\in [0,\LU]$, where $\AU$ denotes
the maximum loss of any event corresponding to peril group $g$ and $\LU$
denotes the maximum yearly ceded loss of peril group $g$. Since $\attach_g\geq
\AL$, the constraint on the attachment probability is automatically satisfied.

For our branch \& bound experiments, we use different data sets. Since we are
interested in the performance of QBB relative to a classical implementation at
different problem sizes, we use a simple model to generate events for each of
the perils at a given problem size (given by the number of perils). For each
trial year $\trial\in\{0,\dots,999\}$, each event loss $\loss_{\trial,\event}$,
where $\event\in\{0,\dots,49\}$ denotes the event index in each year, is drawn
from a Pareto distribution. Since most insurance losses are small and frequent
and the largest losses are rare, the heavy tailed nature of this distribution
effectively models insurance losses. 

Specifically, for each peril group $\group\in\{0,\dots,15\}$ and pair
$(\trial,\event)$, we draw a number $u$ uniformly at random from $[0,1)$ and
set
\begin{equation}
  \loss_{\trial,\event}^{(g)}=s^{(g)} \cdot (1-u)^{-\frac12}\,,
\end{equation}
where $s>1$ is a scaling factor to model the different loss scales of the
different perils. Informed by the scales of the expected losses of the original
data set, we use $s=1.2$ in our tree size experiments.

\subsection{Bounds}
\label{sec:branch-bound-bounds}
In this section, we describe different bounding operators for our simplified
reinsurance optimization problem. In addition to our inexpensive
application-specific bounds, a classical B\&B algorithm could leverage more
expensive bounds, e.g., based on linear programming, which would be too
expensive to implement on a quantum computer~\citep{dalzell2023end}, but could
result in a more efficient classical solution.

%%%%%%%%%%%%%%%%%%%%%%%%%%%%%%%%%%%%%%%%%%%%%%%%%%%%%%%%%%%%%%%%%%%%%%%%%%%%%%%%
\vspace{-5mm}
\subsubsection*{Recursive profit bound:}
As a first step, we derive a recursive bound (similar to the bound used by
\citep{hartwig1984recursive} in the context of Ising spin-glass models) that is
conceptually simple, but computationally expensive due to its recursive nature.
Given fixed layers $\layer_1,\dots,\layer_i$ and variable layers
$\layer_{i+1},\dots,\layer_n$, we can write the net profit $\profit$ as
\begin{equation}\label{eq:profit-bound}
  \profit(\layer_1,\dots,\layer_n)
    = \sum_{j=1}^n \profit(\layer_j)
    \leq \left(\sum_{j\leq i} \profit(\layer_j)\right)
      + \profit_{\max}(\layer_{i+1},\dots,\layer_n)\,,
\end{equation}
where $\profit_{\max}(\layer_{i+1},\dots,\layer_n)$ is the maximum net profit
due to peril groups $i+1,\dots,n$ subject to 
\begin{equation}\label{eq:sp-risk-constraint}
  K(\layer_1, \dots, \layer_i, \layer_{i+1}, \dots, \layer_n) \leq K_{\max}\,.
\end{equation}
We compute $\pi_{\max}$ using our B\&B method, which, in turn, recursively
instantiates B\&B on lower-dimensional sub-problems of
Eq.~\ref{eq:profit-bound} subject to Eq.~\ref{eq:sp-risk-constraint}.

Specifically, solving for $\layer_1,\dots,\layer_n$ starts with a branching
step on $\layer_1$, followed by a bound computation
$\pi(\layer_1)+\pi_{\max}(\layer_2,\dots,\layer_n)$. Determining
$\pi_{\max}(\layer_2,\dots,\layer_n)$ again starts with a branching step and
the bound computation involves $\pi_{\max}(\layer_3,\dots,\layer_n)$, and so
on. Computing the bound at a given node thus involves solving an array of
sub-problems.

Such an array of sub-problems has to be solved at every node of the branch \&
bound tree since all sub-problems---specifically the risk constraint
(Eq.~\ref{eq:sp-risk-constraint})---depend on $\layer_1, \dots, \layer_i$. As
a result, this B\&B approach is costly, mainly due to its recursive nature.

%%%%%%%%%%%%%%%%%%%%%%%%%%%%%%%%%%%%%%%%%%%%%%%%%%%%%%%%%%%%%%%%%%%%%%%%%%%%%%%%
\vspace{-5mm}
\subsubsection*{Iterative profit bound:}
To arrive at a less expensive bounding operator, we now derive an iterative
version of the recursive profit bound. The main advantage of this new bound is
that we have to solve only one optimization problem for each group of perils
$g$ ranging from 1 to $n$. As a result, the solver may first compute auxiliary
solutions $\pi_{\max}(\layer_i,\dots,\layer_n)$ for $i$ ranging from $i=n$ down
to $i=2$, using B\&B with previously computed
$\pi_{\max}(\layer_j,\dots,\layer_n)$ for all $j>i$ as profit upper bounds.
Once all sub-problems have been solved, the final B\&B solves the
full-dimensional problem and returns optimal layers $\layer_1,\dots,\layer_n$.

To derive this bound, we change the risk constraint
(Eq.~\ref{eq:sp-risk-constraint}) such that it becomes independent of previous
branches taken and we fix the branching order globally. We thus redefine
$\pi_{\max}$ to be the maximum profit due to perils (or peril groups)
$i+1,\dots,n$ subject to 
\begin{equation}
  K(\layer_1^*, \dots, \layer_i^*, \layer_{i+1}, \dots, \layer_n) \leq K_{\max}\,,
\end{equation}
where we use the minimum risk solution for the first $i$ towers, denoted by
$\layer_1^*, \dots, \layer_i^*$. Here we set the attachment and limit of
$\layer_g^*$ to $\attach^{\mathrm{min}}_{g}$ and $\limit^{\mathrm{max}}_{g}$.
Note that the constraint on the attachment probability $(P_{\mathrm{attach}})$
makes this minimum risk solution non-trivial, i.e., not all risk is ceded.

While this bound is weaker than the recursive bound, it is more practical to
use the iterative bound in QBB, because it is less expensive to evaluate. In
particular, we do not have to solve a superposition of optimization problems to
compute the bound when implementing the reflection operators of
QBB~\citep{montanaro2020quantum} on a quantum device. In
section~\ref{sec:classical-bb}, we will explore numerically whether or not the
iterative bound is sufficiently strong for B\&B to be effective.

%%%%%%%%%%%%%%%%%%%%%%%%%%%%%%%%%%%%%%%%%%%%%%%%%%%%%%%%%%%%%%%%%%%%%%%%%%%%%%%%%%%%%%%%%
\vspace{-5mm}
\subsubsection*{Risk bound:}
To leverage the risk constraint (Eq.~\ref{eq:risk-constraint}) for pruning
nodes, we also compute a lower bound on the achievable risk. If this lower
bound violates the risk constraint at a given node, we may prune the subtree
rooted at that node.

Since our profit bound is iterative, our solver solves a $k$-peril (group)
problem by first solving $k-1$ lower-dimensional sub-problems, starting with
$g=1$ peril groups, where only the $k$th peril group is considered, up to
$g=k-1$, where all peril groups are considered except the first. For each
sub-problem with $g>1$ peril groups, the solver uses the solutions
$\pi_{\max}(\layer_{n-j+1},\dots,\layer_n)$ of previously solved sub-problems
with $j\in\{1,\dots ,g-1\}$ peril groups to evaluate the profit bound
(Eq.~\ref{eq:profit-bound}) at different depths of the B\&B tree. In addition
to the profit bound, our solver uses a risk lower bound for pruning nodes when
solving these sub-problems and the final problem.

Given a partial (or complete) assignment $(\layer_{k-j+1},\dots,\layer_m)$ with
$m\leq k$, our solver computes the risk lower bound
\begin{equation}\label{eq:risk-bound}
  K_{j,m}^{(LB)}(\layer_{k-j+1}, \dots, \layer_m)
    :=K(\layer_1^\star,\dots,\layer_{k-j}^\star, \layer_{k-j+1},
        \dots, \layer_m, \layer_{m+1}^\star, \dots, \layer_k^\star),
\end{equation}
where $\layer_i^\star$ denotes a minimum risk layer for peril group $i$, and
checks if $K_{j,m}^{(LB)}(\layer_{n-j+1}, \dots, \layer_m) > K_{\max}$. If
true, then $K(\layer_1,\dots,\layer_k)>K_{\max}$ for all
$(\layer_1,\dots,\layer_k)$ that are feasible extensions of the partial
assignment $(\layer_{k-j+1}, \dots, \layer_m)$, allowing the node to be pruned.

%%%%%%%%%%%%%%%%%%%%%%%%%%%%%%%%%%%%%%%%%%%%%%%%%%%%%%%%%%%%%%%%%%%%%%%%%%%%%%%%
\vspace{-3mm}
\subsection{Implementation and results}\label{sec:classical-bb}
We have implemented a B\&B solver leveraging the iterative profit bound and the
risk bound introduced in section~\ref{sec:branch-bound-bounds}, as well as data
compression based on the attachment probability constraint
(Eq.~\ref{eq:attach_prob}), as described in
section~\ref{sec:data-compression-generation}. We decided to use the C++
programming language to enable high performance. To reduce overheads in memory
and runtime due to the data structure representing the branch \& bound tree,
the solver does not store the tree explicitly, but implements the depth-first
tree traversal by only storing and updating a single tree node. These aspects
are important since we expect the branch \& bound trees to be large.

While the discussion of the iterative bound considers fixed layers $\layer_g$
($g\leq i$) and variable layers $\layer_j$ ($j > i$), the actual implementation
is more nuanced. In particular, we cannot afford to exhaustively enumerate all
possible $(\attach,\limit)$-pairs, set the attachment and limit of the current
layer, and then check the bounds. Instead, we perform two interleaved binary
searches for the attachment and limit of the current layer. After each step of
the binary search, we compute an updated lower bound on the total risk
(considering all perils, i.e., assuming minimum-risk solutions for unassigned
peril layers) and check if the node (and the entire subtree rooted at the node)
can be eliminated. Once we arrive at a leaf of the local binary search subtree,
we recompute the profit upper bound as well as the risk lower bound to infer if
we may eliminate the current node, which now corresponds to a
completely-specified layer (attachment and limit have a definite value).

The number of binary search steps $b$ is a solver parameter. Because the number
of leaf nodes of the entire B\&B tree (i.e., the number of possible solutions)
is $2^{2n\cdot b}$, it is currently intractable to consider large values for
$n$ and $b$. In our experiments, we use up to $n=16, b=1$ and $n=5, b=3$. A
stronger bound, different synthetic data sets, or additional constraints (e.g.,
peril-specific risk bounds) would have to be added to allow the solver to find
the (provably) optimal solution for significantly larger instances within
minutes.

%%%%%%%%%%%%%%%%%%%%%%%%%%%%%%%%%%%%%%%%%%%%%%%%%%%%%%%%%%%%%%%%%%%%%%%%%%%%%%%%
\vspace{-5mm}
\subsubsection*{Bound strength:}
To analyze how effective our proposed bound is, we generate and solve 10
instances for each pair $(n,b)$, where $b\in\{1,2,3\}$ and
$n\in\{2,\dots,16/b\}$, and we report the tree size with a trivial bound (i.e.,
without pruning nodes) and the tree size actually observed with the cost and
risk bounds from the previous section.

\begin{figure}[t]
  \begin{center}
    \includegraphics[width=.9\textwidth]{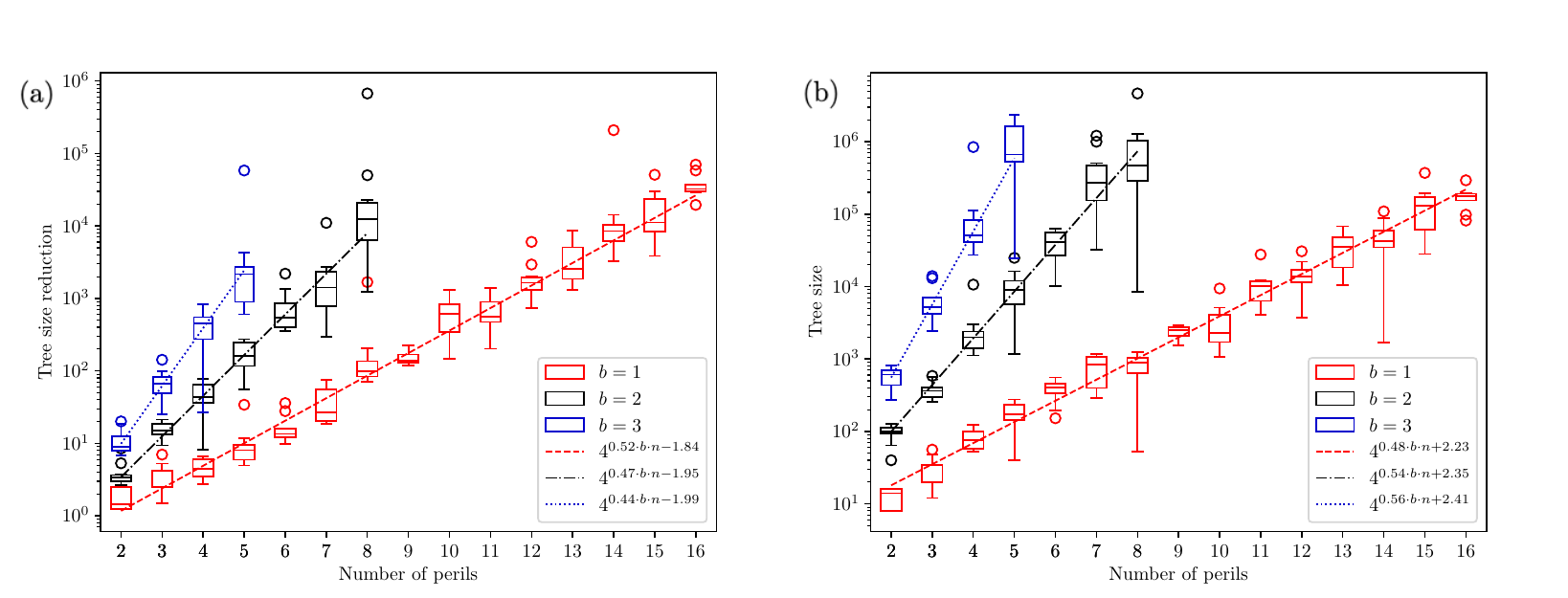}
    \captionsetup{width=.8\linewidth}
    \caption{
      \textbf{(a)} Tree size reduction due to the iterative bound (relative to
        exhaustively enumerating all B\&B tree nodes) as a function of the
        number of peril (groups) $n$ and the number of binary search steps $b$
        for limit and attachment.
      \textbf{(b)} Branch and bound tree sizes as a function of the number of
        peril (groups) $n$ and the number of binary search steps $b$ for limit
        and attachment when using the iterative bound.}
    \label{fig:tree-sizes}
  \end{center}
\end{figure}

We find that the iterative bound effectively eliminates a significant fraction
of the tree; see Fig.~\ref{fig:tree-sizes}(a), which shows the reduction
of the tree size, i.e., the size of the exhaustive search tree divided by the
tree size observed when our bound is used for node elimination. An exhaustive
search would have to consider $2^{2nb}=4^{nb}$ assignments (for each peril, we
do $b$ steps of binary search for both attachment and limit). In contrast, the
B\&B solver visits between $4^{0.48bn+2.23}$ and $4^{0.56bn+2.41}$ nodes, see
Fig.~\ref{fig:tree-sizes}(b) for a similar plot depicting tree sizes. The
bounds thus effectively reduce the number of variables by almost a factor of
$2\times$.

%%%%%%%%%%%%%%%%%%%%%%%%%%%%%%%%%%%%%%%%%%%%%%%%%%%%%%%%%%%%%%%%%%%%%%%%%%%%%%%%
\vspace{-5mm}
\subsection{Quantum computing model and practicality of QBB}
\label{sec:qcmodel}

In this section, we describe in more detail the model of quantum computing that
we consider for making estimates on the practicality of QBB. It is currently
unclear which hardware technology will be best suited to build a fault-tolerant
quantum computer, and which error correcting code will perform best. We focus
on one of the most prominent road maps and outline other approaches that might
improve the speed of future fault-tolerant quantum computers in
section~\ref{sec:discussion}.

We focus on superconducting qubits with a planar surface code for
fault-tolerance. Quantum operations cannot be performed arbitrarily fast: Their
fastest speed for a physical two-qubit gate, referred to as quantum speed limit
of that gate, depends on the physical interaction strength between the two
qubits \citep{PhysRevResearch.5.043194}. Gate times for superconducting qubits
are 12-34$ns$ for the two different two-qubit gates \citep{googleDataSheet}.
The bottlenecks are the 500$ns$ + 160$ns$ for measurement and reset time
\citep{google2023suppressing}. In contrast, trapped ion based quantum computers
require 11-46$\mu s$ for a measurement \citep{crain2019high,todaro2021state}.
This significant difference in measurement times motivates our choice to focus
on superconducting qubits.

For our estimates, we are interested in logical gate times. These depend not
just on the speed of the physical gates and measurements, but also on the
quantum error correction code used to achieve fault tolerance. We assume that
the surface code is used, as it is currently regarded by the community as one
of the most practical error correction schemes for a two-dimensional
architecture with a high error threshold \citep{Fowler2012SurfaceCode}.
\cite{Babbush2021QuadraticSpeedup} estimate that a Toffoli gate will take
around 170$\mu s$ using a distance-30 surface code, assuming that a single
round of surface code including decoding time will be around 1$\mu s$ for
superconducting qubits with a physical error rate of $10^{-3}$
\citep{gidney2019efficient}. Toffoli gate factories can be parallelized when
more qubits are available. A single Toffoli requires at least 10 qubit-seconds
of space-time volume under these assumptions
\citep{Babbush2021QuadraticSpeedup}.

%%%%%%%%%%%%%%%%%%%%%%%%%%%%%%%%%%%%%%%%%%%%%%%%%%%%%%%%%%%%%%%%%%%%%%%%%%%%%%%%
\vspace{-5mm}
\subsubsection*{Number of Toffoli gates per oracle call:}
Given this quantum computing model and the (near\mbox{-})quadratic speedup of
QBB \citep{montanaro2020quantum}, we perform a high-level assessment of
feasibility using the formulas in Eqs.~\ref{eq:t_classical} and
\ref{eq:t_quantum}, as done by \cite{hoefler2023disentangling}. In particular,
we will need to estimate the ratio $t_q/t_c$.

Assuming that we have $10^6$ physical qubits with gate error rate of $10^{-3}$
and require 10 qubit-seconds per Toffoli \citep{Babbush2021QuadraticSpeedup},
we get $10^5$ Toffoli gates per second. For our comparison, we consider
fp16-operations, which may offer sufficient precision for some practical
scenarios. Moreover, we use an optimistic estimation of the quantum resources
required to implement fp16-operations by considering a much simpler $16$-bit
integer addition, which requires $15$ AND gates~\citep{gidney2018halving},
allowing us to estimate that the quantum computer can perform $6666$ operations
per second.

For the classical comparison, we consider a single CPU core running at 3 GHz,
with 2 AVX512 instruction per cycle. This results in $1.9 \cdot 10^{11}$ fp16
FMA (fused multiply add) instructions per second. Note that the quantum 16-bit
integer adder is a significantly simpler operation than the full fp16 FMA
implemented by the classical CPU.

In this quantum-optimistic scenario, the relative speed advantage for 16-bit
floating point operations of a single classical core is $t_q/t_c \approx 2.9
\cdot 10^{7}$. To reach our goal of achieving $T_Q \leq T_C$, we thus need
$\sqrt{N}\geq 2.9 \cdot 10^{7}$, which corresponds to very large tree sizes $N
\gtrsim 8 \cdot 10^{14}$. In addition, we want the calculation to complete in
less than $10^6$ seconds (11.6 days), which is possible only if each call to
the bounding operator requires at most $230$ additions.

Evaluating our bounds in Eqs.\ref{eq:profit-bound} and ~\ref{eq:risk-bound}
using fewer than $230$ additions is clearly infeasible: Even just a single
subtraction for each of the 394,067 events, which is needed to compute the
retained loss for a given attachment, is well above this number. At this point,
we have not taken into account the cost of additional operations needed to
compute the bound (controlled additions, data loading, percentile calculations,
etc.) or the overheads involved in implementing QBB. Moreover, as noted above,
we have only considered a single classical CPU core, instead of allowing for
GPUs or even special-purpose (classical) hardware.
This shows that significant advances in quantum hardware, error correction, and
algorithms are needed to make a future quantum advantage via QBB---or any
quantum algorithm with a quadratic
speedup---feasible~\citep{Babbush2021QuadraticSpeedup,hoefler2023disentangling,haner2024solving}.

%%%%%%%%%%%%%%%%%%%%%%%%%%%%%%%%%%%%%%%%%%%%%%%%%%%%%%%%%%%%%%%%%%%%%%%%%%%%%%%%
%%%%%%%%%%%%%%%%%%%%%%%%%%%%%%%%%%%%%%%%%%%%%%%%%%%%%%%%%%%%%%%%%%%%%%%%%%%%%%%%
\section{Discussion \& Conclusion}
\label{sec:discussion}

% putting things together
In this work, we present and implement advanced computational techniques to
enhance catastrophe excess-of-loss reinsurance contract optimization in
practical settings. We tackle a challenging optimization problem involving
attachment points, limits, covered perils, and optional reinstatement clauses,
all while aiming to maximize expected net profit under internal risk management
and regulatory constraints.

We thoroughly explore the problem formulation and open a path forward for
practitioners through two complementary approaches. Firstly, we introduced a
local search optimizer based on simulated annealing, a heuristic capable of
handling real-world constraints effectively. Secondly, we applied a branch \&
bound approach, an exact algorithm, to a simplified version of the problem.
This approach sheds light on the future potential of quantum-accelerated branch
\& bound algorithms. 

%%%%%%%%%%%%%%%%%%%%%%%%%%%%%%%%%%%%%%%%%%%%%%%%%%%%%%%%%%%%%%%%%%%%%%%%%%%%%%%%
%\vspace{-5mm}
\subsection{Markov-Chain Optimization and Simulated Annealing}
% intuition
Our goal in applying a local search optimizer was to numerically identify
optimal reinsurance structures in a practical context. In particular, our aim
was to conceive an architecture that potentially supports a wide set of
contract features, objective functions, and constraints as dictated by business
needs and regulatory requirements.

% solution
An approach based on simulated annealing is a promising candidate for these
requirements given its versatility across multiple fields of research. Our
Markov chain, defined by the neighborhood structure in
section~\ref{sec:simulated-annealing}, can explore contracts with a rich set of
features found in catastrophe excess-of-loss contracts. And, when implemented
with rapid objective function evaluation in mind, this heuristic optimization
approach remains feasible even for data sets of considerable size (50,000
years, 42MM events).

% caveats
The objective function (expected net profit) and risk constraints (such as
TVaR, AEP and OEP) reported on here represent one specific use case. However,
we want to point out that the approach does not directly leverage these choices
beyond the fact that they are all computed from the (distribution of) ceded
losses. I.e., the method can be adapted to an alternative selection with
limited adjustments and parameter tuning. 

% future looking
In a concrete business setting, the setup may need to be adjusted to support
the specific reinsurance features and risk constraints relevant to the company.
While this study focuses on excess-of-loss layers, it could equally be used to
consider proportional reinsurance, partial placements, alternative
reinstatement terms or inuring contracts. Furthermore, integration with other
business parts could allow for a more holistic risk assessment that considers
correlated exposure when evaluating constraint thresholds.

%%%%%%%%%%%%%%%%%%%%%%%%%%%%%%%%%%%%%%%%%%%%%%%%%%%%%%%%%%%%%%%%%%%%%%%%%%%%%%%%
%\vspace{-5mm}
\subsection{Quantum Branch \& Bound}
For the quantum branch \& bound approach, the iterative bound developed in
section~\ref{sec:branch-bound-bounds} is both inexpensive and effective, and
the B\&B tree sizes grow rapidly. These are necessary conditions for QBB to be
useful in practice, as discussed in the introduction.
In addition, to determine whether QBB might one day achieve a practical quantum
advantage for reinsurance optimization, a more in-depth exploration of the best
classical algorithms for such problems would be required. If there exist
variations of the problem such that the iterative bound performs well compared
to all other classical algorithms, future advances in quantum hardware and
software may enable a quantum advantage via QBB.

For example, \cite{gidney2024magic} propose a potentially cheaper method to
implement T-gates but their construction is not easy to simulate and might
require experiments on real quantum computers to decide if it works well. They
mention that the effect on the overall runtimes of quantum algorithms is hard
to predict. As a first-order effect, T-gates become cheaper but a second-order
effect might be that the space-time tradeoff of different subroutines might
change, i.e., current adders might get exchanged by faster and T-hungrier
carry-lookahead adders.
Other hardware and error correction proposals might outperform the models used
in this paper. One example would be using a cat qubit encoding
\citep{chamberland2022building}.
% For more technical details, we refer to the AWS blog post by
% \cite{aws2021catqubit}.

Data loading presents another challenge for a quantum branch \& bound approach:
At every node of the B\&B tree, the entire classical data set is loaded to
compute the profit bound. This makes the evaluation of the bound very costly.
Instead of loading this data set, one could envision producing these synthetic
data on the quantum computer on-the-fly, allowing the bound to be evaluated
using fewer gates. However, we note that if (for example) a closed-form
expression for the data is available, one would also have to check that there
are no other (classical) methods to solve this problem more efficiently.

Finally, the hardness of the optimization problem depends crucially on the
constraints in the problem definition. In order to identify areas where QBB
might be beneficial, one can investigate if there is a set of constraints for
which our classical B\&B algorithm beats all other classical alternatives such
as (mixed) integer programming and interior point methods.

%%%%%%%%%%%%%%%%%%%%%%%%%%%%%%%%%%%%%%%%%%%%%%%%%%%%%%%%%%%%%%%%%%%%%%%%%%%%%%%%
%\vspace{-5mm}
\subsection*{Concluding remarks}

In summary, this work provides a practical blueprint on how to leverage two
complementary search methods for insurance professionals to optimize
reinsurance contracts. One approach, heuristic and suited to current real-world
applications, provides a computationally viable framework considering realistic
constraints. The other one, exact but currently intractable, outlines a
forward-looking perspective on how reinsurance optimization could be further
advanced through the future adoption of quantum computing technologies.

%%%%%%%%%%%%%%%%%%%%%%%%%%%%%%%%%%%%%%%%%%%%%%%%%%%%%%%%%%%%%%%%%%%%%%%%%%%%%%%%
%%%%%%%%%%%%%%%%%%%%%%%%%%%%%%%%%%%%%%%%%%%%%%%%%%%%%%%%%%%%%%%%%%%%%%%%%%%%%%%%

\subsection*{Acknowledgements}

We thank Pete Martin, Jean Marty, Antoine Pietri and Charles-David Teboul for
their support.

\bibliographystyle{apalike}
\bibliography{references}

\end{document}

%% file: figures/layered.tex
\begin{tikzpicture}
    \node at (-0.65,3.3) {\small{(a)}};
    \node[rotate=90] at (-0.6,1.5) {\small{coverage}};
    \draw[black,thin,-stealth] (-0.3,0) -> (-0.3,3.5);
    \draw[black,thin] (-0.3,0) -> (2.5,0);
    \fill[blue!5!white] (0,0.02) rectangle (2,3);
    \fill[green!70!black!5!white] (0.02,0.6) rectangle (1.98,1.4);
    \node[green!70!black] at (1,0.95) {\small{Layer 1}};
    \fill[teal!5!white] (0.02,1.4) rectangle (1.98,2.99);
    \node[teal] at (1,2.1) {\small{Layer 2}};
    \node at (2.5,0.6) {\small{3MM}};
    \node at (2.5,1.4) {\small{8MM}};
    \node at (2.5,2.99) {\small{16MM}};
    \node at (0.99,0.29) {\small{Retention}};
    \draw[black,thin] (0.02,1.4) -> (1.98,1.4);
    \draw[black,thin] (0.02,0.6) -> (1.98,0.6);
    \draw[black,thin] (0.02,0.6) -> (0.02,2.99);
    \draw[black,thin] (1.98,0.6) -> (1.98,2.99);
    \draw[black,thin] (0.02,2.99) -> (1.98,2.99);
\end{tikzpicture}
\hspace{3mm}
\begin{tikzpicture}
    \node at (-0.65,3.3) {\small{(b)}};
    \node[rotate=90] at (-0.6,1.5) {\small{coverage}};
    \draw[black,thin,-stealth] (-0.3,0) -> (-0.3,3.5);
    \draw[black,thin] (-0.3,0) -> (2.5,0);
    \fill[blue!5!white] (0,0.02) rectangle (2,3);
    \fill[green!70!black!5!white] (0.02,0.6) rectangle (1.98,1.4);
    \fill[teal!5!white] (0.02,1.4) rectangle (1.98,2.99);
    \fill[teal!5!white] (0.02,1.4) rectangle (1.98,2.99);
    
    \draw[black,thin] (0.02,1.4) -> (1.98,1.4);
    \draw[black,thin] (0.02,0.6) -> (1.98,0.6);
    \draw[black,thin] (0.02,0.6) -> (0.02,2.99);
    \draw[black,thin] (1.98,0.6) -> (1.98,2.99);
    \draw[black,thin] (0.02,2.99) -> (1.98,2.99);
    
    \draw[red,dashed, fill=red!15!, opacity=0.5] (0,0.02) rectangle (2,0.85);
    \draw[red,dashed] (-0.3,0.85) -- (2.1,0.85);
    \node[red] at (0.99,1.1) {\small{loss}};
    \node[red] at (2.6,0.85) {\small{5MM}};
\end{tikzpicture}
\hspace{3mm}
\begin{tikzpicture}
    \node at (-0.65,3.3) {\small{(c)}};
    \node[rotate=90] at (-0.65,1.5) {\small{coverage}};
    \draw[black,thin,-stealth] (-0.3,0) -> (-0.3,3.5);
    \draw[black,thin] (-0.3,0) -> (2.5,0);
    \fill[blue!5!white] (0,0.02) rectangle (2,3);
    \fill[green!70!black!5!white] (0.02,0.6) rectangle (1.98,1.4);
    \fill[teal!5!white] (0.02,1.4) rectangle (1.98,2.99);
    \fill[teal!5!white] (0.02,1.4) rectangle (1.98,2.99);
    
    \draw[black,thin] (0.02,1.4) -> (1.98,1.4);
    \draw[black,thin] (0.02,0.6) -> (1.98,0.6);
    \draw[black,thin] (0.02,0.6) -> (0.02,2.99);
    \draw[black,thin] (1.98,0.6) -> (1.98,2.99);
    \draw[black,thin] (0.02,2.99) -> (1.98,2.99);
    
    \draw[red,dashed, fill=red!15!, opacity=0.5] (0,0.02) rectangle (2,1.7);
    \draw[red,dashed] (-0.3,1.7) -- (2.1,1.7);
    \node[red] at (0.99,1.95) {\small{loss}};
    \node[red] at (2.65,1.7) {\small{10MM}};
\end{tikzpicture}

%% file: figures/shifted.tex
\hspace{-1mm}

\begin{tikzpicture}
    \node at (-1.15,3.3) {\small{(a)}};
    \node[rotate=90] at (-1.15,1.5) {\small{coverage}};
    \draw[black,thin,-stealth] (-0.8,0) -> (-0.8,3.5);
    \draw[black,thin,-stealth] (-0.8,0) -> (4.9,0);
 
    % Connecting Tower
    %\fill[fill=teal,fill=teal!5!white,rotate around={-35.8:(1.81,2.99)}] (1.81,2.005) rectangle (3.16,2.99);
    %\fill[fill=green!70!black!5!white,rotate around={-35.8:(1.81,1.85)}] (1.81,1.16) rectangle (3.16,1.8);
    %\fill[fill=blue!5!white,rotate around={-35.8:(1.81,1)}] (1.81,0.5) rectangle (2.6,1);
    
    % First Tower
    \fill[blue,fill=blue!5!white] (0.35,0.01) rectangle (2.08,2.99);
    \fill[fill=green!70!black!5!white] (0.35,1) rectangle (2.08,1.8);
    \fill[teal,fill=teal!5!white] (0.35,1.82) rectangle (2.08,2.99);
    \draw[black,thin] (0.35,1) -> (0.35,1.8);
    \draw[black,thin] (0.35,1.8) -> (0.35,2.99);
    \draw[black,thin] (0.35,1) -> (2.08,1);
    \draw[black,thin] (0.35,1.8) -> (2.08,1.8);
    \draw[black,thin] (0.35,2.99) -> (2.08,2.99);
    \node at (-0.2,1) {\small{3MM}};
    \node at (-0.2,1.8) {\small{8MM}};
    \node at (-0.2,2.99) {\small{16MM}};

    % Second Tower
    \fill[blue,fill=blue!5!white] (1.35,0.01) rectangle (3.81,0.6);
    \fill[green!70!black,fill=green!70!black!5!white] (2.08,0.6) rectangle (3.81,1.4);
    \fill[teal,fill=teal!5!white] (2.08,1.41) rectangle (3.81,2.61);
    \draw[black,thin] (2.08,0.6) -> (3.81,0.6);
    \draw[black,thin] (2.08,1.4) -> (3.81,1.4);
    \draw[black,thin] (2.08,2.62) -> (3.81,2.62);
    \draw[black,thin] (3.81,0.6) -> (3.81,1.4);
    \draw[black,thin] (3.81,1.4) -> (3.81,2.62);    
    \node at (4.4,0.6) {\small{1MM}};
    \node at (4.4,1.4) {\small{6MM}};
    \node at (4.4,2.62) {\small{14MM}};

    \draw[blue!15!white,thin,dashed] (2.08,0) -> (2.08,2.6);
    \draw[black,thin] (2.08,1) -> (2.08,0.6);
    \draw[black,thin] (2.08,1.8) -> (2.08,1.4);
    \draw[black,thin] (2.08,2.99) -> (2.08,2.62);
   
    \node at (1.22,.35) {\small{Retention}};
    \node[green!70!black] at (1.2,1.35) {\small{Layer 1}};
    \node[teal] at (1.2,2.3) {\small{Layer 2}};
    \node at (1.23,-.3) {\small{$s_{1}$}};
    \node at (3.09,-.3) {\small{$s_{2}$}};
\end{tikzpicture}
\hspace{2mm}
\begin{tikzpicture}

    \node at (-1.15,3.3) {\small{(b)}};
    \node[rotate=90] at (-1.15,1.5) {\small{coverage}};
    \draw[black,thin,-stealth] (-0.8,0) -> (-0.8,3.5);
    \draw[black,thin,-stealth] (-0.8,0) -> (4.9,0);
 
    % Connecting Tower
    %\fill[fill=teal,fill=teal!5!white,rotate around={-35.8:(1.81,2.99)}] (1.81,2.005) rectangle (3.16,2.99);
    %\fill[fill=green!70!black!5!white,rotate around={-35.8:(1.81,1.85)}] (1.81,1.16) rectangle (3.16,1.8);
    %\fill[fill=blue!5!white,rotate around={-35.8:(1.81,1)}] (1.81,0.5) rectangle (2.6,1);
    
    % First Tower
    \fill[blue,fill=blue!5!white] (0.35,0.01) rectangle (2.08,2.99);
    \fill[fill=green!70!black!5!white] (0.35,1) rectangle (2.08,1.8);
    \fill[teal,fill=teal!5!white] (0.35,1.82) rectangle (2.08,2.99);
    \draw[black,thin] (0.35,1) -> (0.35,1.8);
    \draw[black,thin] (0.35,1.8) -> (0.35,2.99);
    \draw[black,thin] (0.35,1) -> (2.08,1);
    \draw[black,thin] (0.35,1.8) -> (2.08,1.8);
    \draw[black,thin] (0.35,2.99) -> (2.08,2.99);

    % Second Tower
    \fill[blue,fill=blue!5!white] (1.35,0.01) rectangle (3.81,0.6);
    \fill[green!70!black,fill=green!70!black!5!white] (2.07,0.6) rectangle (3.81,1.4);
    \fill[teal,fill=teal!5!white] (2.08,1.41) rectangle (3.81,2.61);
    \draw[black,thin] (2.08,0.6) -> (3.81,0.6);
    \draw[black,thin] (2.08,1.4) -> (3.81,1.4);
    \draw[black,thin] (2.08,2.62) -> (3.81,2.62);
    \draw[black,thin] (3.81,0.6) -> (3.81,1.4);
    \draw[black,thin] (3.81,1.4) -> (3.81,2.62);    
   
    \draw[blue!15!white,thin,dashed] (2.08,1.6) -> (2.08,2.6);
    \draw[black,thin] (2.08,1) -> (2.08,0.6);
    \draw[black,thin] (2.08,1.8) -> (2.08,1.4);
    \draw[black,thin] (2.08,2.99) -> (2.08,2.62);

    \node at (1.23,-.3) {\small{$s_{1}$}};
    \node at (3.09,-.3) {\small{$s_{2}$}};
    
    % Top Limit
    % Side Lables
    \draw[orange,dashed, fill=orange!15!, opacity=0.5] (2.08,0.01) rectangle (3.81,1.6);
    \draw[orange,dashed] (3.82,1.6) -- (4.0,1.6);
    \node[orange] at (3,1.9) {\small{loss 2}};
    \node[orange] at (4.45,1.6) {\small{7MM}};

    \draw[red,dashed, fill=red!15!, opacity=0.5] (0.35,0.01) rectangle (2.08,1.8);
    \draw[red,dashed] (0.2,1.8) -- (0.35,1.8);
    \node[red] at (1.2,2.1) {\small{loss 1}};
    \node[red] at (-0.25,1.8) {\small{8MM}};
\end{tikzpicture}

\hspace{2mm}

%% file: figures/moves.tex
\begin{tikzpicture}

  \begin{scope}[shift={(0,0.2)}]
    %\node at (-0.65,3.3) {\small{(a)}};
    \fill[teal!5!white] (0.02,1.4) rectangle (1.98,2.6);
    \node[teal] at (0.95,1.95) {\small{Layer 2}};
    \fill[green!70!black!5!white] (0,0.6) rectangle (1.98,1.4);
    \node[green!70!black] at (0.95,0.95) {\small{Layer 1}};
    \fill[blue!5!white] (0,0) rectangle (2,0.6);
    %\node at (1,0.28) {\small{Retention}};
    % horizontals
    \draw[black,thin] (0,2.6) -> (2,2.6);
    \draw[black,thin] (0,1.4) -> (2,1.4);
    \draw[black,thin] (0,0.6) -> (2,0.6);
    \draw[blue!15!white,thin] (0,0) -> (2,0);
    % sides
    \draw[blue!15!white,thin] (0,0) -> (0,0.6);
    \draw[blue!15!white,thin] (2,0) -> (2,0.6);
    \draw[black,thin] (0,0.6) -> (0,2.6);
    \draw[black,thin] (2,0.6) -> (2,2.6);    
  \end{scope}

  \begin{scope}[shift={(-5,2.5)}]
    \node at (-0.05,2.4) {\small{(i)}};
    \fill[teal!5!white] (0.35,1.4) rectangle (3.67,2.6);
    \node[teal] at (1.15,1.95) {\small{Layer 2}};
    \fill[green!70!black!5!white] (0.35,0.6) rectangle (3.67,1.4);
    \node[green!70!black] at (1.15,0.95) {\small{Layer 1}};
    \fill[blue!5!white] (0.35,0) rectangle (3.67,0.6);
    %\node at (1,0.28) {\small{Retention}};
    % horizontals
    \draw[black,thin] (0.35,2.6) -> (3.67,2.6);
    \draw[black,thin] (0.35,1.4) -> (3.67,1.4);
    \draw[black,thin] (0.35,0.6) -> (3.67,0.6);
    \draw[blue!15!white,thin] (0.35,0) -> (3.67,0);
    % sides
    \draw[blue!15!white,thin] (0.35,0) -> (0.35,0.6);
    \draw[blue!15!white,thin] (3.67,0) -> (3.67,0.6);
    \draw[black,thin] (0.35,0.6) -> (0.35,2.6);
    \draw[black,thin] (3.67,0.6) -> (3.67,2.6);    
    \draw[blue!15!white,thin,dashed] (2.01,0) -> (2.01,2.6);
    
    \node at (1.18,-.3) {\small{$s_{1}$}};
    \node at (2.82,-.3) {\small{$s_{2}$}};
  \end{scope}

  \begin{scope}[shift={(0,3.5)}]
    \node at (-0.4,1.4) {\small{(ii)}};
    \fill[green!70!black!5!white] (0,0.6) rectangle (1.98,1.4);
    \node[green!70!black] at (0.95,0.95) {\small{Layer 1}};
    \fill[blue!5!white] (0,0) rectangle (2,0.6);
    %\node at (1,0.28) {\small{Retention}};
    % horizontals
    %\draw[black!40!white,thin,dashed] (0,2.6) -> (2,2.6);
    \draw[black,thin] (0,1.4) -> (2,1.4);
    \draw[black,thin] (0,0.6) -> (2,0.6);
    \draw[blue!15!white,thin] (0,0) -> (2,0);
    % sides
    \draw[blue!15!white,thin] (0,0) -> (0,0.6);
    \draw[blue!15!white,thin] (2,0) -> (2,0.6);
    \draw[black,thin] (0,0.6) -> (0,1.4);
    \draw[black,thin] (2,0.6) -> (2,1.4);    
    \draw[black!40!white,thin,dashed] (0,1.4) -> (0,2);
    \draw[black!40!white,thin,dashed] (2,1.4) -> (2,2);    
  \end{scope}

  \begin{scope}[shift={(3.5,2.5)}]
    \node at (-0.4,2.4) {\small{(iii)}};
    \fill[orange!5!white] (0.02,2) rectangle (1.98,2.6);
    \node[orange] at (1,2.3) {\small{Layer 3}};
    \fill[teal!5!white] (0.02,1.4) rectangle (1.98,2);
    \node[teal] at (1,1.7) {\small{Layer 2}};
    \fill[green!70!black!5!white] (0,0.6) rectangle (1.98,1.4);
    \node[green!70!black] at (1,0.95) {\small{Layer 1}};
    \fill[blue!5!white] (0,0) rectangle (2,0.6);
    %\node at (1,0.28) {\small{Retention}};
    % horizontals
    \draw[black,thin] (0,2.6) -> (2,2.6);
    \draw[black!70!white,thin,dashed] (0,2) -> (2,2);
    \draw[black,thin] (0,1.4) -> (2,1.4);
    \draw[black,thin] (0,0.6) -> (2,0.6);
    \draw[blue!15!white,thin] (0,0) -> (2,0);
    % sides
    \draw[blue!15!white,thin] (0,0) -> (0,0.6);
    \draw[blue!15!white,thin] (2,0) -> (2,0.6);
    \draw[black,thin] (0,0.6) -> (0,2.6);
    \draw[black,thin] (2,0.6) -> (2,2.6);    
  \end{scope}

  \begin{scope}[shift={(-5,-2)}]
    \node at (-0.05,2.8) {\small{(vi)}};
    % Connecting Tower
    %\fill[fill=teal,fill=teal!5!white,rotate around={-35.8:(1.81,2.99)}] (1.81,2.005) rectangle (3.16,2.99);
    %\fill[fill=green!70!black!5!white,rotate around={-35.8:(1.81,1.85)}] (1.81,1.16) rectangle (3.16,1.8);
    %\fill[fill=blue!5!white,rotate around={-35.8:(1.81,1)}] (1.81,0.5) rectangle (2.6,1);
    
    % First Tower
    \fill[blue,fill=blue!5!white] (0.35,0.01) rectangle (2.01,2.99);
    \fill[fill=green!70!black!5!white] (0.35,1) rectangle (2.01,1.8);
    \fill[teal,fill=teal!5!white] (0.35,1.82) rectangle (2.01,2.99);
    \draw[black,thin] (0.35,1) -> (0.35,2);
    \draw[black,thin] (0.35,1.8) -> (0.35,2.99);
    \draw[black,thin] (0.35,1) -> (2.01,1);
    \draw[black,thin] (0.35,1.8) -> (2.01,1.8);
    \draw[black,thin] (0.35,2.99) -> (2.01,2.99);

    % Second Tower
    \fill[blue,fill=blue!5!white] (2.00,0.01) rectangle (3.67,0.6);
    \fill[green!70!black,fill=green!70!black!5!white] (2.01,0.6) rectangle (3.67,1.4);
    \fill[teal,fill=teal!5!white] (2.01,1.41) rectangle (3.67,2.61);
    \draw[black,thin] (2.01,0.6) -> (3.67,0.6);
    \draw[black,thin] (2.01,1.4) -> (3.67,1.4);
    \draw[black,thin] (2.01,2.62) -> (3.67,2.62);
    \draw[black,thin] (3.67,0.6) -> (3.67,1.4);
    \draw[black,thin] (3.67,1.4) -> (3.67,2.62);    
    \draw[blue!15!white,thin] (0.35,0) -> (3.67,0);
    \draw[blue!15!white,thin] (0.35,0) -> (0.35,1);
    \draw[blue!15!white,thin] (3.67,0) -> (3.67,0.6);
    \draw[blue!15!white,thin,dashed] (2.01,0) -> (2.01,2.6);
    %\draw[blue!15!white,thin,dashed] (2.345,0) -> (2.345,0.6);

    \draw[black,thin] (2.01,1) -> (2.01,0.6);
    \draw[black,thin] (2.01,1.8) -> (2.01,1.4);
    \draw[black,thin] (2.01,2.99) -> (2.01,2.62);
   
    %\node at (1.3,.35) {\small{Retention}};
    \node[green!70!black] at (1.15,1.35) {\small{Layer 1}};
    \node[teal] at (1.15,2.3) {\small{Layer 2}};
    \node at (1.18,-.3) {\small{$s_{1}$}};
    \node at (2.82,-.3) {\small{$s_{2}$}};
  \end{scope}

  \begin{scope}[shift={(3.5,-1.5)}]
    \node at (-0.4,2.4) {\small{(iv)}};
    \fill[teal!5!white] (0.02,1.8) rectangle (1.98,2.6);
    \node[teal] at (0.95,2.15) {\small{Layer 2}};
    \fill[green!70!black!5!white] (0,0.6) rectangle (1.98,1.8);
    \node[green!70!black] at (0.95,0.95) {\small{Layer 1}};
    \fill[blue!5!white] (0,0) rectangle (2,0.6);
    %\node at (1,0.28) {\small{Retention}};
    % horizontals
    \draw[black,thin] (0,2.6) -> (2,2.6);
    \draw[black,thin] (0,1.8) -> (2,1.8);
    \draw[black!40!white,thin,dashed] (0,1.4) -> (2,1.4);
    \draw[black,thin] (0,0.6) -> (2,0.6);
    \draw[blue!15!white,thin] (0,0) -> (2,0);
    % sides
    \draw[blue!15!white,thin] (0,0) -> (0,0.6);
    \draw[blue!15!white,thin] (2,0) -> (2,0.6);
    \draw[black,thin] (0,0.6) -> (0,2.6);
    \draw[black,thin] (2,0.6) -> (2,2.6);    
  \end{scope}

  \begin{scope}[shift={(0,-2.8)}]
    \node at (-0.4,2.1) {\small{(v)}};
    \fill[teal!5!white] (0.02,1.1) rectangle (1.98,2.3);
    \node[teal] at (0.95,1.65) {\small{Layer 2}};
    \fill[green!70!black!5!white] (0,0.6) rectangle (1.98,1.09);
    \node[green!70!black] at (0.95,0.8) {\small{Layer 1}};
    \fill[blue!5!white] (0,0) rectangle (2,0.6);
    %\node at (1,0.28) {\small{Retention}};
    % horizontals
    \draw[black,thin] (0,2.3) -> (2,2.3);
    \draw[black,thin] (0,1.1) -> (2,1.1);
    \draw[black!40!white,thin,dashed] (0,1.4) -> (2,1.4);
    \draw[black,thin] (0,0.6) -> (2,0.6);
    \draw[blue!15!white,thin] (0,0) -> (2,0);
    % sides
    \draw[blue!15!white,thin] (0,0) -> (0,0.6);
    \draw[blue!15!white,thin] (2,0) -> (2,0.6);
    \draw[black,thin] (0,0.6) -> (0,2.3);
    \draw[black,thin] (2,0.6) -> (2,2.3);    
  \end{scope}

  \draw[black,thin,stealth-stealth] (1,2.9) -> (1,3.35);
  \draw[black,thin,stealth-stealth] (1,0.05) -> (1,-0.4);
  \draw[black,thin,stealth-stealth] (2.4,2.3) -> (3,2.7);
  \draw[black,thin,stealth-stealth] (2.4,0.3) -> (3,-0.1);
  \draw[black,thin,stealth-stealth] (-0.4,2.3) -> (-1,2.7);
  \draw[black,thin,stealth-stealth] (-3.0,2.05) -> (-3.0,1.35);

\end{tikzpicture}